\documentclass[12pt]{article}
\usepackage{amssymb}
\usepackage{eufrak}
\usepackage{amsmath}
\usepackage{xr}

\title{Differential Galois Theory of Linear Difference Equations}
\author{  Charlotte Hardouin\footnote{IWR, Im Neuenheimer Feld 368, 69120 Heidelberg, Germany, charlotte.hardouin@iwr.uni-heidelberg.de.},  Michael F.
Singer\footnote{North Carolina State University, Department of
Mathematics, Box 8205, Raleigh, North Carolina 27695-8205, USA,
singer@math.ncsu.edu. This material is based upon work supported by
the National Science Foundation under Grant No. CCF-0634123.}}

\begin{document}
\date{December 24, 2007}
\def\QED{\hbox{\hskip 1pt \vrule width4pt height 6pt depth 1.5pt \hskip 1pt}}
\newcounter{defcount}[section]
\setlength{\parskip}{1ex}
\newtheorem{thm}{Theorem}[section]
\newtheorem{lem}[thm]{Lemma}
\newtheorem{cor}[thm]{Corollary}
\newtheorem{prop}[thm]{Proposition}
\newtheorem{defin}[thm]{Definition}
\newtheorem{defins}[thm]{Definitions}
\newtheorem{remark}[thm]{Remark}
\newtheorem{remarks}[thm]{Remarks}
\newtheorem{ex}[thm]{Example}
\newtheorem{exs}[thm]{Examples}
\def\GL{{\rm GL}}
\def\SL{{\rm SL}}
\def\Sp{{\rm Sp}}
\def\sl{{\rm sl}}
\def\sp{{\rm sp}}
\def\gl{{\rm gl}}
\def\PSL{{\rm PSL}}
\def\SO{{\rm SO}}
\def\Gal{{\rm Aut}}
\def\Sym{{\rm Sym}}
\def\tildea{\tilde a}
\def\tildeb{\tilde b}
\def\tildeg{\tilde g}
\def\tildee{\tilde e}
\def\tildeA{\tilde A}
\def\calG{{\cal G}}
\def\calH{{\cal H}}
\def\calT{{\cal T}}
\def\calC{{\cal C}}
\def\calP{{\cal P}}
\def\calS{{\cal S}}
\def\calM{{\cal M}}
\def\calF{{\cal F}}
\def\calY{{\cal Y}}
\def\calO{{\cal O}}
\def\curve{{\rm {\bf C}}}
\def\P1{{\rm {\bf P}}^1}
\def\Ad{{\rm Ad}}
\def\ad{{\rm ad}}
\def\Aut{{\rm Aut}}
\def\Int{{\rm Int}}
\def\CX{{\mathbb C}}
\def\QX{{\mathbb Q}}
\def\NX{{\mathbb N}}
\def\ZX{{\mathbb Z}}
\def\DX{{\mathbb D}}
\def\AX{{\mathbb A}}
\def\sd{{\sigma\der}}
\def\SDP{{\Sigma\Delta\Pi}}
\def\SD{{\Sigma\Delta}}
\def\der{{\partial}}
\def\dd{{\partial}}
\def\disp{{\rm disp}}
\def\pdisp{{\rm pdisp}}
\def\kbar{{\overline{k}}}
\def\ubar{{\tilde{u}}}
\def\zbar{{\tilde{z}}}
\def\semi{\hbox{${\vrule height 5.2pt depth .3pt}\kern -1.7pt\times $ }}

\newenvironment{prf}[1]{\trivlist
\item[\hskip \labelsep{\bf
#1.\hspace*{.3em}}]}{~\hspace{\fill}~$\square$\endtrivlist}
\newenvironment{proof}{\begin{prf}{Proof}}{\end{prf}}
 \def\square{\QED}
 \newenvironment{proofofthm}{\begin{prf}{Proof of Theorem~\ref{ext}}}{\end{prf}}
 \def\square{\QED}
\newenvironment{sketchproof}{\begin{prf}{Sketch of Proof}}{\end{prf}}

\def\calu{{\cal U}}
\def\si{\sigma}
\def\cee{\CX}
\def\HX{{\mathbb H}}
\def\GX{{\mathbb G}}
\def\Krd{{\rm Kr.dim}}
\def\dcl{{\rm dcl}}
\maketitle

\maketitle
\begin{abstract}{We present a Galois theory of difference equations designed to measure the {\em differential} dependencies among solutions of linear difference equations.  With this we are able to reprove H\"older's theorem that the Gamma function satisfies no polynomial differential equation  and are able to give general results that imply, for example,  that no differential relationship holds among solutions of certain classes of $q$-hypergeometric functions.}
\end{abstract}

\newpage
\tableofcontents

\newpage
\section{Introduction}\label{intro} In 1887, Otto H\"older \cite{hoelder} proved that the Gamma function $\Gamma(x)$ satisfies no differential polynomial equation, that is, there is no nonzero polynomial $P(x, y, y', \ldots )$ such that $P(x, \Gamma(x), \Gamma'(x), \ldots ) = 0$.  This result has been reproved and generalized over  the years by many researchers (for example, \cite{bankkaufman, carmichael, hausdorff, moore, ostrowski25,
rosenlicht_rank}; see also \cite{rubel_trans}).  Most recently, Hardouin (\cite{hardouin06} and \cite{hardouin07}) (and subsequently van der Put, see the appendix to \cite{hardouin07}) proved and generalized this result using the Galois theory of difference equations (as developed in particular in \cite{PuSi}).  This Galois theory associates a linear algebraic group to a linear difference equation and, using  properties of linear algebraic groups, Hardouin was able  to derive her results.\\[0.1in]
 In this paper we develop a general Galois theory of difference and differential equations where the Galois groups are linear {\em differential} groups, that is groups of matrices whose entries lie in a differential field and satisfy a set of polynomial differential equations. This general theory encompasses the usual Galois theory of linear differential equations \cite{PuSi2003}, the Galois theory of linear difference equations \cite{PuSi} and the Galois theory of parameterized differential equations \cite{CaSi}\footnote{In  \cite{andre_galois}, Andr\'e develops a Galois theory that encompases the difference and differential Galois theories and considers the differential Galois theory as a limiting case of the difference Galois theory.  Our theory does not have that feature while Andr\'e's theory does not address the questions studied in this paper. Another approach to describing analytic properties of solutions of difference equations involving pseudogroups has been announced by Umemura in \cite{um06}.}. We will develop this theory in its full generality in Section~\ref{appendix}.  We will use this theory in a restricted setting, that is, when we wish to analyze the differential behavior of a solution of a linear difference equation and describe this restricted  theory and the tools we need for subsequent sections in Section~\ref{galoistheory}. \\[0.1in]
 In Section~\ref{relations}, we describe possible differential relations among solutions of linear difference equations.  The key idea is that the form of possible differential relations among solutions of linear difference equations is determined by the form of the differential equations defining the associated Galois group.   We begin by considering first order equations.  We prove in our setting a general result  which implies the following result ({\em cf.} Corollary~\ref{indefsumcor} below).  This result (and its  $q$-analogue) already appears in  Hardouin's work (\cite{hardouin07}, Prop. 2.7).
  \begin{quotation}\noindent {\em Let $\CX(x)$ be the field of rational functions over the complex numbers and $\calF$ the field of $1$-periodic functions meromorphic on the complex plane. Let $a_1(x), \ldots , a_n(x) \in \CX(x)$ and let $z_1(x), \ldots , z_n(x)$ be functions, meromorphic on $\CX$ (resp. $\CX^*$) such that
\[z_i(x+1) -z_i(x) = a_i(x), \  i = 1, \ldots ,n.\]
The functions $z_1(x), \ldots ,z_n(x)$ are differentially dependent over $\calF(x)$   if and only if there exists a nonzero homogeneous linear differential polynomial $L(Y_1, \ldots ,Y_n) $ with coefficients in $\CX$ such that $L(a_1(x), \ldots , a_n(x)) = g(x+1) - g(x)$.
 }\end{quotation}
 We also give a similar result for difference equations of the form $z(x+1)=a(x)z(x)$ and   $q$-difference versions (as does Hardouin in \cite{hardouin06}).  These results can be considered as an analogue of the Kolchin-Ostrowski Theorem \cite{kolchin_ostrowski} which characterizes the possible algebraic relations among solutions of first order differential equations. This latter result follows (using the Picard-Vessiot Theroy) from a description of the algebraic subgroups of products of one dimensional linear algebraic groups. In our setting these results follow from a description of the {\em differential} algebraic subgroups of products of one dimensional linear algebraic groups in the same general way once the machinery of our Galois theory  is established.  The theorem of H\"older follows from these results.  Continuing with first order equations, we use facts about solvable differential subgroups of $\GL_2$ to reprove and generalize ({\em cf.} Propositions~\ref{Ishiprop}, \ref{shift_diff_alg}, and \ref{q_diff_alg} below) the following result of Ishizaki \cite{ishizaki}.
 \begin{quotation}\noindent{\em If $a(x), b(x) \in \CX(x)$ and $z(x) \notin \CX(x)$ satisfies $z(qx) = a(x) z(x) + b(x), \ |q|\neq 1$ and is meromorphic on $\CX$, then $z(x)$ is not differentially algebraic over $\calG(x)$, where $\calG$ is the field of $q$-periodic functions meromorphic on $\CX^*$.}\end{quotation}

 We then turn to higher order equations.  Using facts about differential subgroups of simple algebraic groups, we can characterize differential relationships among solutions of difference equations whose difference Galois group is a simple, noncommutative, algebraic group ({\em cf.} Proposition~\ref{higher_order_prop} and Corollary~\ref{higher_order_cor}). Using this and Roques's \cite{roques} computation of the difference Galois groups of the $q$-hpergeometric equations,  we can show, for example ({\em cf.} Example~\ref{q-diffex}),
 \begin{quotation} \noindent {\em Let $y_1(x), y_2(x)$ be linearly independent solutions of the hypergeometric equation
 \[y(q^2 x) - \frac{2ax -2}{a^2x-1}y(qx) + \frac{x-1}{a^2x-q^2}y(x) = 0\]
 where $a\notin q^\ZX$ and $a^2 \in q^\ZX$ and $|q| \neq 1$. Then $y_1(x), y_2(x), y_1(qx)$ are differentially independent over $\calG(x)$, where $\calG$ is the field of $q$-periodic functions meromorphic on $\CX^*$.}
 \end{quotation}
 In Section~\ref{inverse}, we consider a special case of the  problem of determining which linear differential algebraic groups occur as Galois groups. In Section~\ref{param}, we consider certain parameterized families $Y(qx,t) = A(x,t)Y(x,t)$ of difference equations and show that the associated connection matrix (\cite{etingof}) is independent of $t$ if and only if the Galois group we associate with this equation ({\em a priori} a linear differential algebraic group) is conjugate to a linear algebraic group. The paper ends with an appendix containing two subsections,  In Section~\ref{diffappendix} we gather together some facts about  rational solutions of difference equations that are used throughout the preceding sections and, as mentioned above, in Section~\ref{appendix} we   present our general Galois theory of linear difference and differential equations.\\[0.1in]
 We wish to thank Carsten Schneider for references to the literature on finding rational solutions of difference equations
 as well as Daniel Bertrand for
 useful comments and  advice.

\section{Galois Theory}\label{galoistheory}

In this section we will give the basic definitions and results needed in Sections \ref{relations} - \ref{param}.  These results follow from a  more general approach to the Galois theory of differential and difference equations that we  present in Section~\ref{appendix}, where complete proofs are also given. The parenthetical references indicate the relevant general statements and results from the appendix.

\begin{defin} (Definition~\ref{def6.1}) A {\em $\sd$-ring} is a commutative ring $R$ with unit together with an automorphism $\sigma$ and a derivation $\der$ satisfying $\sigma(\der(r)) = \der(\sigma(r)) \ \forall r \in R$. A $\sd$-field is defined similarly\footnote{All fields considered in this paper are of characteristic $0$.}.
\end{defin}

\begin{exs} {\em $1. \ \CX(x), \sigma(x) = x+1, \der(x) = \frac{d}{dx},\\
2. \  \CX(x), \sigma(x) = qx, \ (q\in \CX \backslash\{0\}), \der = x\frac{d}{dx},\\
3. \ \CX(x,t), \sigma(x) = x+1,
\sigma(t) = t, \der = \frac{\partial}{\partial t}$.}
\end{exs}

For $k$ a $\sigma\der$-field we shall consider difference equations of the form
\begin{eqnarray}\label{msequation1}\sigma(Y) = A Y, \ \ A \in \GL_n(k)\end{eqnarray}

\begin{defin} (Definition~\ref{defPV}) A {\em $\sigma\der$-Picard-Vessiot ring ($\sigma\der$-PV-ring)} over $k$ for equation (\ref{msequation1}) is a $\sigma\der$-ring $R$ containing $k$ satisfying: \begin{enumerate}
\item $R$ is a simple $\sigma\der$-ring, {\em i.e.,} $R$ has no ideals, other than $(0)$ and $R$,  that are invariant under $\sigma$ and $\der$
\item There exists a matrix $Z\in \GL_n(R)$ such that $\sigma(Z) = AZ$.
\item $R$ is generated as a $\der$-ring over $k$ by the entries of $Z$ and $1/\det(Z)$, {\em i.e., $R = k\{Z,1/\det(Z)\}_\der$}
\end{enumerate}
\end{defin}

Note that when $\der$ is identically zero, this corresponds to the usual definition of a Picard-Vessiot extension for a difference equation.  To prove existence and uniqueness of Picard-Vessiot extensions, one needs to assume that the field of $\sigma$-invariant elements of $k$ is algebraically closed.  In the case of $\sigma\der$-PV extensions, $k^\sigma = \{ c \in k \ | \ \sigma(c) = c\}$ is a differential field with derivation $\der$ and we need to assume that this field is {\em differentially closed} (see Section 9.1 of \cite{CaSi} for the definition and references.) With this assumption, we have

\begin{prop} (Propositions~\ref{nonewconstants} and~\ref{unique})  Let $k$ be a $\sd$-field with $k^\sigma$  a differentially closed field.
 There exists a $\sd$-PV ring for (\ref{msequation1}) and it is unique up to $\sd$-$k$-isomorphism. Furthermore, $R^\sigma = k^\sigma$.
\end{prop}

\begin{defin}  The {\em $\sd$-Galois group} $\Gal_\sd(R/k)$ of  the $\sd$-PV ring $R$ (or of  (\ref{msequation1})) is \[\Gal_\sd(R/k)=\{\phi \ | \ \mbox{ $\phi$ is a $\sd$-$k$-automorphism of $R$}\} \ .\]
\end{defin}

As in the usual theory of linear difference equations, once one has selected a fundamental solution matrix of (\ref{msequation1}) in $R$, this group may be identified with elements of $\GL_n(k^\sigma).$  The next result states that both  $\Gal_\sd(R/k)$ and $R$ have additional structure.

\begin{thm}\label{torsor} (Propositions~\ref{galoisgp} and~\ref{torsorprop})  Let $k$ be a $\sd$-field and assume that $k^\sigma$ is a differentially closed field. Let $R = k\{Z,\frac{1}{\det Z}\}_\der, \ \sigma(Z) = AZ$ be a $\sd$-PV extension of $k$.
 \begin{enumerate}
 \item We may identify $\Gal_\sd(R/k)$ with the set of $k^\sigma$-points of a linear $\der$-differential algebraic group $G$ defined over $k^\sigma$.
 \item $R$ is a reduced ring and is the  coordinate ring of a $G$-torsor $V$ defined over $k$. The action of $G$ on $V$ induces an action  of $G(k^\sigma)$ on $R$ that corresponds to the action of $\Gal_\sd(R/k)$ on $R$ under the above identification.
 \end{enumerate}
 \end{thm}

 In the above result we use the terms ``coordinate ring'' and ``$G$-torsor'' in the differential sense (see Sections 4 and 9.4 of \cite{CaSi} for definitions of these terms as well as other definitions, facts and references concerning linear differential algebraic groups.)  In the Appendix, we will  furthermore show (Lemma~\ref{idempotents}) that $R = R_0\oplus \ldots \oplus R_{t-1}$ is the finite direct sum of integral differential $k$-algebras $R_i$ where $\sigma:R_{i{\rm mod} t}\rightarrow R_{i+1{\rm mod} t}$ isomorphically.  In particular, the quotient fields of the $R_i$ all have the same $\der$-differential transcendence degree over $k$ (see \cite{DAAG}, Ch. II.10).  Abusing language, we define this to be the {\em $\der$-differential dimension, $\der{\rm -dim.}_k(R)$ of $R$ over $k$}.  The above theorem implies that the $\der{\rm -dim.}_k(R/k)$ is the same as the differential dimension of (the identity component) of $G$, $\der{\rm -dim.}_C(G)$ (see Proposition~\ref{dimprop}).

Finally, we have the usual Galois correspondence (cf., Theorem 1.29 of \cite{PuSi}). We note that if $R$ is a $\sd$-PV-extension of $k$ and $K$ is the total ring of quotients, then any $\sd$-$k$-automorphism of $K$ must leave $R$ invariant.  Therefore one can identify the group of $\sd$-$k$-automorphisms $\Gal_\sd(K/k)$ of $K$ with the $\sd$-Galois group $\Gal_\sd(R/k)$.

\begin{thm} (Theorem~\ref{fundthm}) Let $k$ and $R$ be as in Theorem~\ref{torsor} and let $K$ be the total ring of quotients of $R$.  Let $\calF$ denote the set of $\sd$-rings $F$ with $k \subset F \subset K$ such that every non-zerodivisor of $F$ is a unit in $F$ and let $\calG$ denote the set of $\der$-differential subgroups of $\Gal_\sd(R/k)$.  There is a bijective correspondence $\alpha: \calF \rightarrow \calG$ given by $\alpha(F) = G(K/F) = \{\phi \in \Gal_\sd(K/k) \ | \ \phi(u) = u \ \forall u \in F\}$.  The map $\beta : \calG \rightarrow \calF$ given by $\beta(H) = \{ u \in K \ | \ \phi(u) = u \ \forall \phi \in H\}$ is the inverse of $\alpha$.
\end{thm}

In particular, an element of $K$ is left fixed by all $\phi$ in $\Gal_\sd(K/k)$ if and only if it is in $k$ and  for a differential subgroup $H$ of $\Gal_\sd(K/k)$ we have $H = \Gal_\sd(K/k)$ if and only if $K^H = k$.

In the next section we will need the following consequences of the above result.  Let $R= k\{Z,\frac{1}{\det Z}\}_\der$ be a $\sd$-PV-ring where $\sigma(Z) = AZ$.  We can consider the $\sigma$-ring $S= k[Z,\frac{1}{\det Z}] \subset R$.  Note that this ring is not necessarily closed under the action of $\der$.  In the next result we show that $S$ is the usual PV ring (as in \cite{PuSi}) for the difference equation $\sigma(Y) = AY$. To avoid confusion we will use the prefix ``$\sigma$-'' to denote objects ({\em e.g.}, $\sigma$-PV-extension, $\sigma$-Galois group) from the Galois theory of difference equations described in \cite{PuSi}.

\begin{prop}\label{sdzariskidense} (Proposition~\ref{propcompare}) Let $k$ and $R$ be as in Theorem~\ref{torsor} and let $S = k[Z,\frac{1}{\det Z}] \subset R$. \begin{enumerate}
\item $S$ is the $\sigma$-PV-extension of $k$ corresponding to $\sigma(Y) = AY$, and
\item $\Gal_\sd(R/k)$ is Zariski dense in the $\sigma$-Galois group $\Gal_\sigma(S/k)$.

\end{enumerate}

\end{prop}

 The following result characterizes those difference equations whose $\sd$-Galois groups are ``constant''.

 \begin{prop} \label{constantgaloisgroup} Let  $k$ and $R $ be as in Theorem~\ref{torsor} and let $C = k^\sd = \{ c \in k \ | \ \sigma(c) = c \mbox{ and } \der(c) = 0\}$.  The $\sd$-Galois group $\Gal_\sd(R/k) \subset \GL_n(k^\sigma)$ is conjugate over $k^\sigma$ to a subgroup of $\GL_n(C)$ if and only if there exists a $B \in \gl_n(k)$ such that
 \[\sigma(B) = ABA^{-1} + \der(A) A^{-1} \ .\]
 In this case, there is a solution $Y=U \in \GL_n(R)$ of the system
 \begin{eqnarray*}
 \sigma(Y) & = & AY\\
 \dd(Y) & = & BY
 \end{eqnarray*}

 \end{prop}

 \begin{proof} Assume that such a $B$ exists.  A calculation shows that $\sigma(\der(Z) - BZ) = A(\der(Z) - BZ)$. Therefore $\der(Z) - BZ = DZ$ for some $D \in \gl_n(k^\sigma)$. We therefore have $\der(Z) - (B+D) Z = 0$. For any $\phi \in \Gal_\sd(R/k)$ we will denote by $[\phi]_Z \in \GL_n(k^\sigma)$ the matrix such that $\phi(Z) = Z[\phi]_Z$.   We claim that $[\phi]_Z \in \GL_n(C)$.  To see this note that $\der(\phi(Z)) = (B+D)Z[\phi]_Z + Z\der([\phi]_Z)$ and $\phi(\der(Z)) = (B+D)Z[\phi]_Z$.  This implies that $Z\der([\phi]_Z) = 0$ so $\der([\phi]_Z)=0$.\\[0.1in]
 Now assume that there exists a $D \in \GL_n(k^\sigma)$ such that $D^{-1} \Gal_\sd(R/k)D \subset \GL_n(C)$.  For $\phi \in \Gal_\sd(R/k)$, let $[\phi]_Z \in \GL_n(k^\sigma)$ once again be the matrix such that $\phi (Z) = Z[\phi]_Z$.  Let $U = ZD$.  For any $\phi \in \Gal_\sd(R/k)$ we have that $\phi(U) = Z[\phi]_ZD = ZD(D^{-1}[\phi]_ZD)$.  Therefore $\phi(U) = U[\phi]_U$ for some $[\phi]_U \in \GL_n(C)$.  This implies that $B = \der(U)U^{-1}$ is left fixed by $\Gal_\sd(R/k)$ and so $B \in \gl_n(k)$.  A calculation shows that $\sigma(B) = \sigma(\der(U)U^{-1}) = ABA^{-1} + \der(A) A^{-1}$.\end{proof}

\section{Differential Relations Among Solutions of Difference Equations}\label{relations}  In this section we shall show how the Galois theory of Section \ref{galoistheory} can be used to give necessary and sufficient conditions for solutions of linear difference equations to satisfy differential polynomial equations.

\subsection{First order equations.}\label{firstorder} The classical Kolchin-Ostrowski Theorem \cite{kolchin_ostrowski} implies that if $k\subset K $ are differential fields with the same constants $C$ and $z_1, \ldots ,z_n \in K$ with $z_i' = a_i \in k$, then $z_1, \ldots , z_n$ are algebraically dependent over $k$ if and only if there exists a homogeneous linear polynomial $L(Y_1, \ldots , Y_n)$ with coefficients in $C$ such that $L(z_1, \ldots , z_n)= f \in k$ or, equivalently, $L(a_1, \ldots , a_n) = f', \ f\in k$.  Kolchin proved this using differential Galois theory and the fact that algebraic subgroups of $(C^n, +)$ are precisely the vector subspaces. For difference equations, the analogy of the indefinite integrals above are indefinite sums, that is, elements $y$ satisfying $\sigma(y) - y = a$ and one has a similar result characterizing algebraic dependence.  The following characterizes {\em differential} dependence among indefinite sums.  Recall that if $k \subset K$ are differential rings, we say that $z_1, \ldots , z_n \in K$ are differentially dependent over $k$ if there is a nonzero differential polynomial $P \in k\{Y_1, \ldots, Y_n\}$ such that $P(z_1, \ldots ,z_n) = 0$ (cf., \cite{DAAG}, Ch. II).

\begin{prop}\label{indefsum} Let $k$ be a $\sd$-field with $k^{\sigma}$ differentially closed and let  $S \supset k$ be a $\sd$-ring such that $S^\sigma = k^\sigma$.  Let $a_1, \ldots , a_n \in k$ and $z_1, \ldots , z_n \in S$ satisfy
\[ \sigma(z_i) - z_i = a_i \ \ i = 1, \ldots, n \ .\]
Then $z_1, \ldots , z_n$ are differentially dependent over $k$ if and only if there exists a  nonzero homogeneous linear differential polynomial $L(Y_1, \ldots , Y_n)$ with coefficients in $k^\sigma$ and an element $f \in k$ such that
\[L(a_1, \ldots , a_n) = \sigma(f) - f.\]

\end{prop}

\begin{proof}  Assuming there exists such an $L$, one sees that $L(z_1, \ldots z_n) - f$ is left fixed by $\sigma$ and so lies in $k^\sigma$.  This yields a relation of differential dependence over $k$ among the $z_i$.\\[0.1in]
Now assume that the $z_i$ are differentially dependent over $k$.  Differentiating the relations $\sigma(z_i) - z_i = a_i$, one sees that the $\der$-differential ring $k\{z_1, \ldots ,z_n\}_\der$ generated by the $z_i$ is a $\sd$-ring so we may assume that $S = k\{z_1, \ldots ,z_n\}_\der$.  Let $M$ be a maximal $\sd$-ideal in $S$ and let $R = S/M$.  $R$ is a simple $\sd$-ring.
Furthermore, $R =k\{Z, \frac{1}{\det Z}\}_\der$ where   $\sigma(Z) = AZ$ and
\[A = \left(\begin{array}{ccccccc}
1&a_1&0&0&\cdots& 0&0\\
0&1&0&0&\cdots&0&0\\
0&0&1&a_2&\cdots&0&0\\
0&0&0&1&\cdots&0&0\\
\vdots&\vdots&\vdots&\vdots&\vdots&\vdots&\vdots\\
0&0&\cdots&\cdots&\cdots&1&a_n\\
0&0&\cdots&\cdots&\cdots&0&1
\end{array}\right) \mbox{ and }
Z = \left(\begin{array}{ccccccc}
1&\overline{z}_1&0&0&\cdots& 0&0\\
0&1&0&0&\cdots&0&0\\
0&0&1&\overline{z}_2&\cdots&0&0\\
0&0&0&1&\cdots&0&0\\
\vdots&\vdots&\vdots&\vdots&\vdots&\vdots&\vdots\\
0&0&\cdots&\cdots&\cdots&1&\overline{z}_n\\
0&0&\cdots&\cdots&\cdots&0&1
\end{array}\right)\]
with $\overline{z}_i$ being the image of $z_i$ in $R$. The $\sd$-Galois group $\Gal_\sd(R/k)$  is a differential subgroup of $(k^\sigma, +)^n$.  By assumption, the differential dimension of $R$ over $k$ is less than $n$, the differential dimension of $(k^\sigma, +)^n$ and so $\Gal_\sd(S/k)$ is a proper differential subgroup of $(k^\sigma, +)^n$. Cassidy has classified these groups (\cite{cassidy1, CaSi}): the differential subgroups of $(k^\sigma, +)^n$ are precisely the zero sets of systems of homogeneous linear differential polynomials over $k^\sigma$.  Therefore there exists a nonzero homogeneous linear differential polynomial $L(Y_1, \ldots , Y_n)$ with coefficients in $k^\sigma$ such that $\Gal_\sd(R/k) \subset \{(c_1, \ldots , c_n) \in (k^\sigma, +)^n \ | \ L(c_1, \ldots , c_n) = 0\}$.

We claim that $L(\overline{z}_1, \ldots , \overline{z}_n) = f \in k$.  To prove this it is enough to show that this element is left fixed by $\Gal_\sd(R/k)$.  Let $\phi \in \Gal_\sd(R/k)$.  We have that $\phi(L(\overline{z}_1,\ldots , \overline{z}_n) = L(\overline{z}_1 + c_1,\ldots , \overline{z}_n+c_n) = L(\overline{z}_1. \ldots , \overline{z}_n)$ so the claim is proved.  Finally we have that $L(a_1, \ldots, a_n) = L(\sigma(\overline{z}_1)- \overline{z}_1, \ldots , \sigma(\overline{z}_n)- \overline{z}_n) = \sigma(f) - f$.
\end{proof}
The above result has the rather artificial assumption that $k^\sigma$ is differentially closed.  Nonetheless, this result can be used to prove results about meromorphic functions. In the following, we denote by $\calF$ the field of $1$-periodic meromorphic functions, that is meromorphic functions $f(x)$ on $\CX$ such that $f(x+1) = f(x)$ and let $\calG$ be the field of $q$-periodic meromorphic functions, that is functions $f(x)$, meromorphic on $\CX^* = \CX\backslash \{0\}$ such that $f(qx) = f(x)$, where $q \in \CX, |q| \neq 1$. In the following we shall speak of homogeneous linear differential polynomials $L(Y_1, \ldots, Y_n)$, that is,  linear forms in the variables $Y_i^{(j)}$. When we substitute elements $a_i$ of a differential field $(k,\dd)$ for the variables $Y_i$, we  will replace $Y_i^{(j)}$ with $\dd^j(a_i)$.  In particular, when we are considering the shift $\sigma(x) = x+1$, we will use the derivation $\dd = \frac{d}{dx}$ and when we consider $q$-difference equations we will use the derivation $\dd = x\frac{d}{dx}$.

\begin{cor}\label{indefsumcor}  Let $a_1(x), \ldots , a_n(x) \in \CX(x)$ and let $z_1(x), \ldots , z_n(x)$ be functions, meromorphic on $\CX$ (resp. $\CX^*$) such that
\[z_i(x+1) -z_i(x) = a_i(x), \ \mbox{(resp. $z_i(qx) - z_i(x) = a_i(x))$ for } i = 1, \ldots ,n.\]
The functions $z_1(x), \ldots ,z_n(x)$ are differentially dependent over $\calF(x)$ (resp. $\calG(x)$)  if and only if there exists a nonzero homogeneous linear differential polynomial $L(Y_1, \ldots ,Y_n) $ with coefficients in $\CX$ such that $L(a_1(x), \ldots , a_n(x)) = g(x+1) - g(x)$  (resp. $L(a_1(x), \ldots , a_n(x)) = g(qx) - g(x))$  for some  $g(x) \in \CX(x)$.

\end{cor}

\begin{proof}We will deal with the case of the shift and apply Proposition~\ref{indefsum} with $\sigma(x) = x+1$ and $\dd = \frac{d}{dx}$.  The $q$-difference case is similar except that we have $\sigma(x) = qx$ and $\dd = x\frac{d}{dx}$. Clearly if $L(a_1(x), \ldots , a_n(x)) = g(x+1) - g(x)$ for some $g(x) \in \CX(x)$, then $L(z_1(x), \ldots , z_n(x)) - g(x) \in \calF$ so the $z_i$ are differentially dependent over $\calF(x)$.\\[0.1in]
Now assume that the $z_i$ are differentially dependent over $\calF(x)$.  Note that $\calF(x)$ is a $\sd$ field with $\sigma(x) = x+1$ and $\der = \frac{d}{dx}$.  Furthermore, the function $f(x) = x$ is not algebraic over $\calF$ since any polynomial equation over $\calF$ satisfied by $x$ would also be satisfied by $x+n$ for all $n \in \ZX$.  Let $T = \calF(x)\{z_1(x), \ldots , z_n(x)\}_\der$. This is a $\sd$-domain and $T^\sigma = \calF$. Let $\calC$ be the differential closure of $\calF$ and define a $\sd$-structure on $R = T\otimes_\calF\calC$ via $\sigma(t\otimes c) = \sigma(t)\otimes c$ and $\der(t\otimes c) = \der(t)\otimes c + t\otimes \der(c)$.  We note that $R^\sigma = \calC$.\\[0.1in]
Letting $k = \calC(x)$ and $R$ as above, we apply
Proposition~\ref{indefsum}.  We can conclude that there exists a
nonzero homogeneous linear differential polynomial $\tilde{L}(Y_1,
\ldots ,Y_n) $ with coefficients in $\calC$ such that
$\tilde{L}(a_1(x), \ldots , a_n(x)) = \tilde{g}(x+1) - \tilde{g}(x)
\mbox{ for some } \tilde{g}(x) \in \calC(x)$. Replace the
coefficients of the variables in $\tilde{L}(Y_1, \ldots , Y_n)$ by
indeterminates to get a differential polynomial $\bar{L}$ with
indeterminate coefficients.  Replace the coefficients of powers of
$x$ in $\tilde{g}$ with indeterminates to get a rational function
$\bar{g}(x)$ with indeterminate coefficients.  We know that there is
a specialization in $\calC$ of the indeterminates so the equation
$\bar{L}(a_1(x), \ldots, a_n(x)) = \bar{g}(x+1)-\bar{g}(x)$ is
satisfied.  By clearing denominators and equating like powers of
$x$, we see that this latter equation is equivalent to a system of
polynomial equations in the indeterminates.  Since $\CX$ is
algebraically closed and this system has a solution in some
extension, it has a solution in $\CX$. Specializing to this
solution, yields an $L$ and $g$ that satisfies the conclusion of the
corollary. \end{proof}
There is also a multiplicative version of the above corollaries. Again, we note that in the shift case $\dd = \frac{d}{dx} $ and for $q$-difference equations, $\dd = x\frac{d}{dx}$.
\begin{cor}\label{hypergeomcor} Let $b_1(x), \ldots , b_n(x) \in \CX(x)$ and let $u_1(x), \ldots , u_n(x)$ be nonzero functions, meromorphic on $\CX$ (resp. $\CX^*$) such that
\[u_i(x+1) = b_i(x)u_i(x), \ \mbox{(resp. $u_i(qx) = b_i(x)u_i(x)$) for } i = 1, \ldots ,n.\]
The functions $u_1(x), \ldots ,u_n(x)$ are differentially dependent over $\calF(x)$ (resp.$\calG(x)$)  if and only if there exists a nonzero homogeneous linear differential polynomial $L(Y_1, \ldots ,Y_n) $ with coefficients in $\CX$ such that $L(\frac{\dd(b_1(x))}{(b_1(x)}, \ldots , \frac{\dd(b_n(x))}{b_n(x)}) = g(x+1) - g(x)$  (resp. $L(\frac{\dd(b_1(x))}{(b_1(x)}, \ldots , \frac{\dd(b_n(x))}{b_n(x)})$  $= g(qx) - g(x))$ for some  $g(x) \in \CX(x)$..
\end{cor}

\begin{proof} Again, we only prove the corollary in the case of the shift. Let $z_i(x) = \frac{u_i'(x)}{u_i(x)}$.  Since the domain $\calF(x)\{u_1, \ldots , u_n\}_\der$ is differentially algebraic over $\calF(x)\{z_1, \ldots ,z_n\}_\der$, standard facts concerning differential transcendence degree imply that the $z_i$ are differentially dependent over $\calF(x)$ if and only if the $u_i$ are differentially dependent over $\calF(x)$.  The $z_i$ satisfy
\[z_i(x+1) - z_i(x) = \frac{b_i'(x)}{b_i(x)}.\]
Corollary~\ref{indefsumcor} implies the conclusion. \end{proof}
We will need the following special case of Th\'eor\`eme 4.12 and Proposition 5.1 of \cite{hardouin06} and shall prove this using the present techniques.
\begin{cor}\label{hypergalois} Let $b(x) \in \CX(x)$.\\[0.1in]
1.~If $u(x)$ is a nonzero function meromorphic on  $\CX$ satisfying $u(x+1) = b(x) u(x)$, then $u(x)$ is differentially algebraic over $\calF(x)$ if and only if $b = c\frac{f(x+1)}{f(x)}$for some $c\in \CX$ and $f(x) \in \CX(x)$.\\[0.1in]
2.~If $u(x)$  is a nonzero function meromorphic on  $\CX^*$ satisfying $u(qx) = b(x) u(x)$, then $u(x)$ is differentially algebraic over $\calG(x)$ if and only if $ b= cx^n\frac{f(qx)}{f(x)}$ for some $a \in \CX, n \in \ZX$ and $f(x) \in \CX(x)$.\end{cor}
\begin{proof} Some standard facts concerning rational solutions of difference equations are contained in Section~\ref{diffappendix}.  We shall freely refer to these in the following proof.\\[0.1in]
 1.~Lemma~\ref{normalform} implies that we may write $b = \tilde{b}\frac{\sigma(\tilde{g})}{\tilde{g}}$ where $\tilde{g} \in \CX(x)$ and  $\tilde{b}$ is standard, that is, the set of zeroes and poles of $\tilde{b}$ contains no numbers that differ by nonzero integers. One sees that the element $w = ug^{-1}$ satisfies $w(x+1) = \tilde{b}w(x)$ and that $w$ is differentially algebraic if and only if $u$ is.  Therefore we shall assume that $b$ is standard and show that it is constant if and only if $u$ is differentially algebraic over $\calF(x)$.   From Corollary~\ref{hypergeomcor}, we see that $u$ is differentially algebraic if and only if there exists a nonzero $L \in \CX[\frac{d}{dx}]$ and such that $L(\frac{b'(x)}{b(x)}) = g(x+1) - g(x)$ for some $g \in \CX(x)$.  Lemma~\ref{trivlem} implies that $\frac{b'(x)}{b(x)} = e(x+1) - e(x)$ for some $e \in \CX(x)$.  If $b$ is nonconstant, then $\frac{b'(x)}{b(x)}$ has a pole.  This implies that  $e$ must have a pole and Lemma~\ref{posdisp} implies that two poles of $\frac{b'(x)}{b(x)}$ differ by a nonzero integer, yielding a contradiction.  If $b$ is constant, we see that $\frac{u'}{u}$ is in $\calF$.\\[0.1in]
2.~ We again can reduce to the case where $b$ is standard.  We will
show that in this case, $u$ is differentially algebraic over
$\calG(x)$ if and only if   $b = cx^n$ for some $c\in \CX$ and $n
\in \ZX$.  From Corollary~\ref{hypergeomcor}, we see that $u$ is
differentially algebraic if and only if there exists a nonzero $L
\in \CX[x\frac{d}{dx}]$  such that $L(\frac{xb'(x)}{b(x)}) = g(qx) -
g(x)$ for some $g \in \CX(x)$.  Lemma~\ref{trivlem} implies that
$\frac{xb'(x)}{b(x)} = e(qx) - e(x)$ for some $e \in \CX(x)$. If $b
\neq cx^n$ then $\frac{xb'(x)}{b(x)}$ will have a pole and we are
led to a contradiction as before.  If $b = cx^n$, then a calculation
shows that $x(\frac{xu'}{u})' \in \calG$.\end{proof}

We note that Praagman \cite{praagman} has shown that any difference equation \linebreak $Y(x+1) = A(x) Y(x), \ A(x) \in \GL_n(\CX(x))$  has a solution meromorphic on $\CX$  and that any difference equation $Y(qx) = AY(x), \ |q|\neq 1$ has a solution meromorphic in $\CX^*$. Therefore  the hypotheses of the above corollaries  are not vacuous.

Corollary~\ref{hypergalois}.1 is easily seen to imply that the Gamma function  $\Gamma(x)$, (where $\Gamma(x+1) = x \Gamma(x)$) satisfies no nonzero polynomial differential equation over $\calF(x)$.

Corollary~\ref{hypergeomcor} appears in  \cite{hardouin06, hardouin07} where two proofs (due to Hardouin and van der Put)  using the usual Galois theory of difference equations are given.   One can also find there a necessary and sufficient condition on the divisors of the $a_i$ (or the $\frac{\dd b_i}{b_i}$) that guarantees the existence of an $L$ and $g$ as in the above corollaries.

We now turn to the following result of Ishizaki \cite{ishizaki} concerning meromorphic solutions of inhomogeneous first order $q$-difference equations.
\begin{prop}\label{Ishiprop} ({\em cf.}~Theorem 1.2 \cite{ishizaki}) If $a(x), b(x) \in \CX(x)$ and $z(x) \notin \CX(x)$ satisfies $z(qx) = a(x) z(x) + b(x)$ and is meromorphic on $\CX$, then $z(x)$ is not differentially algebraic over $\calG(x)$. \end{prop}
We shall prove  more general results (Propositions~\ref{shift_diff_alg} and \ref{q_diff_alg}) that also include first order inhomogeneous difference equations with respect to the shift. Before we do this we need some additional facts concerning differential algebraic groups.\\[0.1in]
 Let $(k, \dd)$ be a differentially closed field ordinary differential field.  We need to consider differential subgroups of
\[\GX = \{ \left(\begin{array}{cc} \alpha & \beta \\ 0 & 1\end{array}\right) \ | \ 0\neq \alpha \in k, \beta \in k\}.\]
We will denote by $\GX_u$ the subgroup
\[\GX_u = \{ \left(\begin{array}{cc} 1 & \beta \\ 0 & 1\end{array}\right) \ | \  \beta \in k\}.\]
\begin{lem}\label{groups} Let $H$ be a differential subgroup of $\GX$.  \begin{enumerate}
\item If $H_u = \GX_u \cap H$ is a proper subgroup of $\GX_u$ then there exists a nonzero linear differential operator $L \in k[\dd]$ such that
\[\left(\begin{array}{cc} 1 & \beta \\ 0 & 1\end{array}\right) \in H_u\]
 if and only if $L(\beta) = 0$.
 \item Assume that  $H_u$ is a proper subgroup of $\GX_u$ that contains more that one element. If
 \[\left(\begin{array}{cc} \alpha & \beta \\ 0 & 1\end{array}\right) \in H\]
then $\dd(\alpha) = 0$.\end{enumerate}
\end{lem}
\begin{proof} 1.~ The group $\GX_u$ is isomorphic to the additive group.  Cassidy \cite{cassidy1} showed that proper differential subgroups of the additive group are of the form described.\\[0.1in]
2. Let $H_u = \{ \left(\begin{array}{cc} 1 & \beta \\ 0 & 1\end{array}\right) \ | \  L(\beta)=0 \}$, where $L$ has order at least one. If $g = \left(\begin{array}{cc} \alpha & \beta \\ 0 & 1\end{array}\right) \in H$ and $h = \left(\begin{array}{cc} 1 & \beta \\ 0 & 1\end{array}\right) \in H_u$, one sees that $ghg^{-1} \in H$ implies that $L(\alpha \beta) = 0$.  Lemma~\ref{linearconst} below implies that $\dd(\alpha) = 0$.\end{proof}
\begin{lem}\label{linearconst} Let $(k,\dd)$ be a differentially closed ordinary differential field and $L \in k[\dd]$ a monic differential operator of order at least one.  Let $\alpha \in k$ and assume that $L(\alpha \beta) = 0$ for all $\beta \in k$ satisfying $L(\beta)=0$.  Then $\dd(\alpha) = 0$. \end{lem}
\begin{proof} Let $y$ and $z$ be differential indeterminates.  A calculation shows that \[L(yz) = zL(y) + nz'y^{(n-1)} + L_{n-2}(z)y^{(n-2)} + \ldots L_1(z)y' +L_0(z)y.\]
where each $L_i$ is an operator of order $n-i$ with coefficients in $k$.  If $\beta_1, \ldots , \beta_n$ are solutions of $L(y) = 0$ independent over the constants then we have \[n\dd(\alpha)\dd^{(n-1)}(\beta_i) + L_{n-2}(\alpha)\dd^{(n-2)}(\beta_i) + \ldots L_1(\alpha)\dd(\beta_i) +L_0(\alpha)\beta_i= 0, \ i=1, \ldots n.\]
Since the wronskian determinant of the $\beta_i$ is nonzero, we must have that $\dd(\alpha) = 0$.\end{proof}
We now turn to a special case of our generalization of Ishizaki's result.   We will consider solutions of an equation of the form $\sigma(z) = az+b, a, b \in k, a\neq 0$ in a $\sd$-PV extension associated with the  equation
\begin{eqnarray}\label{maineqn}
\sigma(Y) = \left(\begin{array}{cc} a & b \\ 0 &1\end{array}\right)Y.
\end{eqnarray}
\begin{prop}\label{constanta} Let $k$ be a $\sd$-field with $k^\sigma$ differentially closed and let $K$ be the total ring of quotients of a $\sd$-PV extension $R$ of $k$ corresponding to equation (\ref{maineqn}) . Let $z \in K$ satisfy
$\sigma(z) = az+b, \ 0\neq a, b \in k$ and assume that $z \notin k$. \begin{enumerate}
\item If the $\sd$-PV group of $\sigma(y) = ay$ over $k$ is $\GL_1(k^\sigma)$, then $z$ is not differentially algebraic over $k$.
\item If $\dd(a) = 0$, then $z$ is differentially algebraic over $k$ if and only if there is a nonzero linear differential operator $L \in k^\sigma[\dd]$ and an element $f \in k$ such that
\[\sigma(f) - af = L(b) \]
\end{enumerate}
\end{prop}
\begin{proof}  Note that by the assumption on $K$, there exists an invertible $u \in K$ such that $\sigma(u) = au$. Furthermore, the  matrix
\[\left(\begin{array}{cc} u & z \\ 0 &1\end{array}\right)\] is a fundamental solution matrix of (\ref{maineqn}).  So we will  assume that $R = k\{u, u^{-1}, z\}_\dd$.  One sees that the   $\sd$-PV group of this equation is a subgroup of $\GX(k^\sigma)$, where $\GX$ is as defined above.\\[0.1in]
1.~Since $z$ is not in $k$, there is a $\sd$-automorphism $\phi$ of $K$ such that $\phi(z) \neq z$.  If $z$ is differentially algebraic over $k$ then $u = \phi(z) - z$ would also be differentially algebraic over $k$. This would imply that the  $\sd$-PV group of $\sigma(u) = au$ would have differential transcendence degree less than $1$ and so be a proper subgroup of $\GL_1(k^\sigma)$.\\[0.1in]
2.~We now assume that $\dd(a) = 0$.  We first note that this implies that the $\sd$-PV group $\tilde{G}$  of $\sigma(y) = ay$   is a subgroup of $\GL_1(C)$ where $C = \{c \in k^\sigma \ | \ \dd(c) = 0\}$.  To see this note that $\sigma(\dd u) = a \dd u$ so $\dd u = d u$ for some $ d \in k^\sigma$. If $\phi \in \tilde{G}$ , then $\phi(u) = \alpha u$ for some $\alpha \in k^\sigma$.  Since $\dd(\phi(u)) = d \phi(u)$, we have $\dd(\alpha) = 0$.  Note that we have also shown that $k[u, u^{-1}]$ is the the $\sd$-PV ring associated with the equation $\sigma(u) = au$. \\[0.1in]
Let $G$ be the $\sd$-PV group of $K$ over $k$ and $H$ be the subgroup of $G$ that leaves the total ring of quotients $F$  of $k[u, u^{-1}]$ elementwise fixed.  Note that $H = G \cap \GX_u$ and so $H$ may be identified with a differential subgroup of $\mathbb{G}_a(k^\sigma)$.   Note that $H$ is a normal subgroup of $G$ and $\tilde{G} = G/H$.  Since the differential transcendence degree  of $K$ over $k$ is equal to the differential dimension of $G$, we have that $z$ is differentially algebraic over $k$ if and only if the differential dimension of $H$ is zero. \\[0.1in]
Assume that $z$ is differentially algebraic over $k$ and therefore that $H$ is a proper differential subgroup of $\GX_u(k^\sigma)$.  Therefore there exists a nonzero monic linear differential operator $\tilde{L} \in k^\sigma[\dd]$ such that
\[ H = \{\left(\begin{array}{cc} 1 & \beta \\ 0 & 1\end{array}\right) \ | \ \tilde{L}(\beta) = 0\}\]
Using the above representation of $H$, one sees that $\tilde{L}(zu^{-1})$ is left fixed by all elements of $H$ and therefore that $\tilde{L}(zu^{-1}) = \tilde{f} \in F$.  Furthermore we have \[\tilde{L}(\frac{b}{au}) = \sigma(\tilde{L}(zu^{-1})) - \tilde{L}(zu^{-1}) = \sigma(\tilde{f}) - \tilde{f}\]
By induction on $i$, one can show that for any $i \geq 0$ we have that $\dd^i(\frac{b}{au}) = L_i(b) (au)^{-1}$ for some  monic $L_i \in k^\sigma[\dd]$. Therefore $\tilde{L}(\frac{b}{au}) = (au)^{-1}L(b)$ for some monic $L \in k^\sigma[\dd]$.  From the equation $(au)^{-1}L(b) = \tilde{L}(\frac{b}{au}) = \sigma(\tilde{f}) - \tilde{f}$ we conclude that $L(b) = \sigma(f) - af$ where $f = u\tilde{f} \in k(u)$.  We will now show how to conclude that there is a $g \in k$ such that $L(b) = \sigma(g) - ag$\\[0.1in]
If $\tilde{G} = G/H$, is finite then we can average  the equation over $\tilde{G}$ to find a \[g  = \frac{1}{|\tilde{G}|}\sum_{\phi \in \tilde{G}} \phi(f) \in k\] such that $L(b) = \sigma(g) - ag.$ \\[0.1in]
Now assume $\tilde{G}$ is infinite.  Since it is a   differential subgroup of $\mathbb{G}_m(C)$, it is really an algebraic group and so must be  connected.  Therefore $k[u]$ is a domain and  $u$ is transcendental over $k$.  Lemma~\ref{descent} yields the conclusion. \\[0.1in]
We now assume that $\dd(a) =0$ and that there exists an $f \in k$ and nonzero $L \in k^\sigma[\dd]$ such that $\sigma(f)-f = L(b)$.  Note that $\sigma(L(z)) - aL(z) = L(b)$.  Therefore there exists an $ 0\neq h \in k^\sigma$ such that
$L(z) - f  = hu$.  Since $\dd(u) = du$ for some $d \in k$, we have
\[ \dd(\frac{1}{h}(L(z) - f) - d(\frac{1}{h}(L(z) - f) = 0.\]
 This shows that $z$ is differentially algebraic over $k$. \end{proof}
 We now turn to the case  $k = k_0(z)$  where  either
 \begin{enumerate}
 \item[(A)] $\dd = \frac{d}{dx}$,  $\sigma(x) = x+1$ and $k_0$ is the field of functions $g$ meromorphic on $\CX$ and satisfying $g(x+1) = g(x)$, or
 \item[(B)] $\dd = x\frac{d}{dx}$,  $\sigma(x) = qx$ where  $q \in \CX$ is not a root of unity and $k_0$ is the field of functions $g$ meromorphic on $\CX\backslash\{0\}$  and satisfying $g(qx) = g(x)$
 \end{enumerate} We shall consider equations of the form $\sigma(y) = ay+b$ with $a \in \CX(x)\backslash\{0\},  b \in \CX(x)$ and discuss when such solutions are differentially algebraic over $k$.  By making a change of variable $y:= fy$ where $f$ is a suitable rational function, we may assume that $a$ is in  standard form ({\em cf.}~Definition~\ref{dispdef}).  This means that if $\alpha$ and $\beta$ are distinct elements from the set of zeroes and poles of $a$, then $\alpha \notin \beta + \ZX$ in case $\sigma(x) = x+1$ or
 $\alpha \notin q^\ZX \beta$ in case $\sigma(x) = qx$.  This change of variable does not affect the property of a solution being differentially algebraic over $k$.\\[0.1in]
\begin{prop}\label{shift_diff_alg} Let $k$ be  as in (A) above and $ a, b \in \CX(x)$ with $a \neq 0$ in standard form.  Let $z$ be a function meromorphic on $\CX$   satisfying $\sigma(z) = a z + b$ and assume that $z \notin k$. Then
\begin{enumerate}
\item If $a \notin \CX$, then $z$ is not differentially algebraic over $k$.
\item If $a \in \CX$, then $z$ is differentially algebraic over $k$ if and only if $b = \sigma(f) - af$ for some $f \in \CX(x)$.
\end{enumerate}
\end{prop}
\begin{proof} Let $u$ be a nonzero solution of $\sigma(y) = ay$ that is meromorphic in $\CX$ (results of \cite{praagman} imply that such a solution exists).  Consider the $\sd$-field $k<z,u>_\dd$, that is, the field generated by $z,u$ and all their derivatives.  This can be seen to be closed under $\sigma$. The ring $k\{u, \frac{1}{u}, z\}_\dd$ can be written as $k\{U, \frac{1}{U}, Z\}_\dd/I$ where $U,Z$ are differential indeterminates and $I$ is a proper $\sd$-ideal.  Let $\kbar_0$ be the differential closure of $k_0$ and extend $\sigma$ and $\partial$ to $\kbar = \kbar_0(x)$ such that $\sigma(f) = f$ for all $f \in \kbar_0$.  The ideal $J = \kbar_0\cdot I$ is again a proper $\sd$-ideal in $\kbar_0\{U,\frac{1}{U}, Z\}_\dd$ and so lies in a maximal $\sd$-ideal $M$. We then have that $R = \kbar_0\{U,\frac{1}{U}, Z\}_\dd/M$ is a $\sd$-PV extension for the equation (\ref{maineqn}). Let $\ubar, \zbar$ be the images of $U,Z$ in $R$ and let $K$ be the total ring of quotients of $R$. Let us assume that $z$ is differentially algebraic over $k$.  We then have that $\zbar$ is differentially algebraic over $\kbar$.\\[0.1in]
Now assume that $a \notin \CX$ and we will derive a contradiction. First, let us show that $\zbar \notin \kbar$.  If $\zbar \in \kbar$ the we may  write
\[\zbar = \frac{\sum_{i=0}^{m-1}a_ix^i}{\sum_{i=0}^{m-1}b_ix^i}\]
where $a_i, b_i \in \kbar_0$.  Clearing denominators of $\sigma(\zbar) = a\zbar + b$ and equating coefficients of powers of $x$, we see that there is a $\CX$-constructible subset $X$ of $\kbar_0^{2m}$ such that $(\tilde{a}_i, \tilde{b}_i) \in X$ if and only if $\tilde{z} = \frac{\sum_{i=0}^{m-1}\tilde{a}_ix^i}{\sum_{i=0}^{m-1}\tilde{b}_ix^i}$.  Since $X$ is nonempty, it has a $\CX$-point and so $\sigma(y) = ay + b$ has a solution $\tilde{z}$ in $\CX(x)$.  The element $z-\tilde{z}$ would be a nonzero solution of $\sigma(y) = ay$ that is differentially algebraic over $k$, contradicting  Corollary~\ref{hypergalois} (or Proposition 5.1 of \cite{hardouin06}).  Therefore, if $a\notin \CX$, $\zbar$ is not in $\kbar$. \\[0.1in]
Note that Corollary~\ref{hypergalois}.1 also implies that $\ubar$ is not differentially algebraic over $\kbar$. Since $\zbar \notin\kbar$ there exists a $\sd$-automorphism $\phi$ of $R$ over $k$ such that $\phi(\zbar) -\zbar \neq 0$.  This element is a solution of $\sigma(y) = ay$ and so must be a $\kbar_0$-multiple of $\ubar$.  Since $\phi(\zbar) -\zbar$ is also differentially algebraic over $k$, we get a contradiction.  This proves 1~above.\\[0.1in]
We now turn to 2.  Assume that $a \in \CX$ and $b = \sigma(f) - af$ for some $f \in \CX(x)$.  Since $e^{z\log a}$ is a solution of $\sigma(y) = ay$, we have that $z = f+ge^{z\log a}$ for some $g \in k_0$.  Therefore $f$ is differentially algebraic over $k$.\\[0.1in]
Now assume that $a \in \CX$ and $z$ is differentially algebraic over $k$. We then have that $\zbar$ is differentially algebraic over $\kbar$. Proposition~\ref{constanta}.3 implies that there is an element $f \in \kbar$ and a nonzero monic linear operator $L\in \kbar_0[\dd]$ such that $\sigma(f) -f = L(b)$. Lemma~\ref{trivlem}.1  implies that $b = \sigma(e) - ae$ for some $e \in \kbar$.  The rational function $e$ has coefficients in $\kbar_0$ and, arguing as above, we can replace this with a rational function $\overline{e}$ with coefficients in $\CX$ such that $b = \sigma(\overline{e}) - a\overline{e}$. \end{proof}

\begin{prop}\label{q_diff_alg} Let $k$ be  as in (B) above and $ a, b \in \CX(x)$ with $a \neq 0$ in standard form.  Let $z$ be a function meromorphic on $\CX$   satisfying $\sigma(z) = a z + b$ and assume that $z \notin k$. Then
\begin{enumerate}
\item If $a\neq cx^n, c\in \CX, n \in \ZX$, then $z$ is not differentially algebraic over $k$.
\item If $a = cx^n$ for some $c \in \CX, n \in \ZX$, then $z$ is differentially algebraic over $k$ if and only if one of the following holds:
\begin{enumerate}
\item $b = \sigma(f) - af + dx^r$ for some $f \in \CX(x), d \in \CX$ when $a = q^r, r \in \ZX$, or
\item $b = \sigma(f) - af $ for some $f \in \CX(x)$, when $a \neq q^r$.
\end{enumerate}
\end{enumerate}
\end{prop}
\begin{proof}  If $a\neq cx^n, c\in \CX, n \in \ZX$, then Corollary~\ref{hypergalois}.2 (or  Th\'eor\`eme 4.12 of \cite{hardouin06}) implies that $\sigma(y) = ay$ has no nonzero differentially algebraic solutions.  One can then argue as in the proof of the first part of  Proposition~\ref{shift_diff_alg} to conclude that $z$ is not differentially algebraic over $k$.\\[0.1in]
Now assume that $a = cx^n$.  Let $u$ be a nonzero solution of $\sigma(y) = cx^ny$.  A calculation shows that $\sigma(\frac{\dd(u)}{u}) - \frac{\dd(u)}{u} = n$. Differentiating again, on sees that $\dd(\frac{\dd(u)}{u})$ must be in $k_0$.  Therefore $u$ is differentially algebraic over $k$. If  $a = q^r$ and $b=\sigma(f) - af + dx^r$ for some $f \in \CX(x), d \in \CX$, we let $L = \dd - r$. One then has that $L(z-f)$ satisfies $\sigma(y) = ay$ and so $L(z-f) = gu$ for some $g \in k_0$.  This implies that $z$ is differentially algebraic over $k$. If $a \neq q^r$ and $b = \sigma(f) - af $ for some $f \in \CX(x)$, then $z-f$ satisfies $\sigma(y) = ay$.  This implies that $z-f$ satisfies $\sigma(y) = ay$ and so $L(z-f) = gu$ for some $g \in k_0$. Again we conclude that $z$ is differentially algebraic over $k$.\\[0.1in]
Finally assume that $a = cx^n$ and that $z$ is differentially algebraic over $k$.  If $n = 0$, we are  in the case covered by Proposition~\ref{constanta}, and so the conclusion follows as in Proposition~\ref{shift_diff_alg}.  We will therefore assume that $n \neq 0$.\\[0.1in]
Proceeding as in the proof of Proposition~\ref{shift_diff_alg}, one constructs a $\sd$-PV extension $R = \kbar\{\ubar, \zbar\}$ of $\kbar$ for the system $\sigma(z) = az+b, \ \sigma(u) = au$, where $\kbar = \kbar_0(x), \kbar_0 =$ the differential closure of $k_0$ and  where $\zbar$ is differentially algebraic over $\kbar$.
Let $G$ be the $\sd$-PV group of $K$ over $\kbar$, where $K$ is the total ring of quotients of $R$ and let $H$ be the $\sd$-PV group of $K$ over $\kbar<\ubar>_\dd$.  Any element of  $H$ leaves $\ubar$ fixed and takes $\zbar$ to $\zbar + \beta \ubar$ for some $\beta \in \kbar^\sigma$. Therefore    $H$ is identified with  a proper subgroup of $\GX_u$.\\[0.1in]
If $\zbar \notin  \kbar<\ubar>_\dd$, then $H$ is nontrivial and so, by Lemma~\ref{groups}.2, if $\phi \in G$, then $\phi(\ubar) = \gamma_\phi \ubar$ where $\gamma_\phi \in \CX$.  This implies that $\dd\ubar/\ubar  = d \in \kbar$. Since $\sigma(\dd(\ubar)) = \dd(\sigma(\ubar))$, we have that $\sigma(d) - d = n$.  Since $d \in \kbar_0(x)$, this can only happen if $n = 0$.  Since we are assuming that $n \neq 0$, we have a contradiction so $\tilde{z} \in \kbar<\ubar>_\dd$\\[0.1in]
We  will now show that there exists an $f \in \kbar$ such that $\sigma(f) - af = b$.  We will then show (as before) that this implies that there exists a $g \in k$ such that $\sigma(g) - ag = b$. \\[0.1in]
Since $\zbar \in \kbar<\ubar>$, we have $K = \kbar<\ubar>$. Note that $\sigma(\ubar) = a\ubar$ and $\sigma(\zbar) = a \zbar + b$. If $G$ is finite, then $g = \frac{1}{|G|}\sum_{\phi\in G}\phi(f)$ satisfies our desired conclusion. If $G$ is infinite, we will first produce the desired $g$ in $F = \kbar<\dd(\ubar)/\ubar>_\dd$ and then show that it must lie in $\kbar$. Since $G$ is infinite,  it is connected (since $G$ is a subgroup of $G_m$, this follows from the Corollary on p.~928 of \cite{cassidy1}).  Therefore, $K$ is a field.   If $\ubar $ is algebraic over $F$, taking traces again yields a desired $g$ in $F$.  If $\ubar$ is not algebraic over $F$, Lemma~\ref{descent} allows us to find a desired $g$ in $F$.  We now claim that if $g \in F$ satisfies $\sigma(g) - ag = b$, then $g \in \kbar$. Note that $w = \dd(\ubar)/\ubar$ satisfies $\sigma(w) - w = n$, so $F$ is a $\sd$-PV extension of $\kbar$ whose $\sd$-Galois group $\tilde{G}$ is a subgroup of $\mathbb{G}_a(\kbar_0)$. \\[0.1in]
If $g \notin \kbar$, there exists a $\phi \in \tilde{G}$ such that $\phi(g) \neq g$.  Therefore $F$ would contain a nonzero solution of $\sigma(y) = ay$ and so contains $\ubar$.  We therefore have  $\kbar<\ubar>_\dd \subset \kbar<w>_\dd$ and that the $\sd$-PV group of $\kbar<\ubar>_\dd $, a subgroup of $\mathbb{G}_m$, is a homomorphic image of a subgroup of $\mathbb{G}_a$.  The image of a unipotent matrix under a differential homomorphism is unipotent ({\em cf.} Proposition 35 of \cite{cassidy1}) so   this image must be the trivial group.  Therefore $\ubar \in \kbar$ and this can only happen if $n=0$, a contradiction. Therefore $g \in \kbar$.\\[0.1in]
Given $g \in \kbar$, we write $g$ as a rational function of $x$ with
coefficients in $\kbar_0$.  The equation $\sigma(g) - ag = b$ shows
that the vector of coefficients is determined by a set of polynomial
equalities and inequalities with coefficients in $\CX$ and, as
before, we can find coefficients in $\CX$.
\end{proof}
We can now give a proof of Ishizaki's result.\\[0.1in]
\noindent {\bf Proof of Proposition~\ref{Ishiprop}.} By replacing $z$ by $zg$ for some $g \in \CX(x)$, we may assume that $a(x)$ is standard and that the hypotheses still hold. Assume that $z(x)$ is differentially algebraic over $k=\calG(x)$.  We shall first show that this implies that $a(x) = cx^n$ for some $c \in \CX$ and $n \in \ZX$ and that $\sigma(y) = ay$ has a solution meromorphic in $\CX$.  \\[0.1in]
Let us first assume that,   $z(x) \in k$. Then the equation $z(qx) =
a(x) z(x) + b(x)$ implies that the coefficients of powers of $x$
appearing in $z$ are determined by  finite sets of polynomial
equalities and inequalities over $\CX$.  Therefore, they will have a
solution in $\CX$ and, using these as coefficients, we have a
solution  $\tilde{z}$ of $\sigma(y) = ay+b$ with $\tilde{z} \in
\CX(x)$.  Therefore $z-\tilde{z}$ is a solution of $\sigma(y) = ay$
that is differentially algebraic over $k$ and is meromorphic in
$\CX$.
Furthermore, Corollary~\ref{hypergalois} implies that $a = cx^n$.\\[0.1in]
Now assume that $z(x) \notin k$.  Proposition~\ref{q_diff_alg} implies that $a =cx^n$ and that $b = \sigma(f) - af +dx^r$ (if $a=q^r$) or $b = \sigma(f) - af$ (if $a\neq q^r$). In the first case, $(\dd-r)(h-f)$ is a nonzero meromorphic solution of $\sigma(y) = ay$ and in the second case, $h-f$ is a nonzero meromorphic solution of this equation. \\[0.1in]
Let $u(x)$ be a function meromorphic in $\CX$ that satisfies $u(qx) = a(x) u(x)$. Since $q^m \rightarrow 0$ or $q^{-m}\rightarrow 0$ as $m \rightarrow \infty$, one sees that $u$ cannot have any poles or zeroes in $\CX^*$, since otherwise, either $u$ would not be meromorphic at $0$ or $u$ would be identically $0$.  Therefore $w = \frac{\dd(u)}{u}$ is a  solution of $\sigma(w) - w= n$ that is meromorphic on $\CX^*$ without any poles.  Differentiating we see that $\dd(w)$ would be a $q$-periodic function without poles.  This function gives a meromorphic function on the elliptic curve $(\CX\backslash\{0\})/q^\ZX$ with no poles and so $\dd(w) = d$ for some $d\in \CX$.  Integrating we see that $w = d \log x$, contradicting the fact that $w$ is meromorphic at $0$.  \hfill\QED

\subsection{Higher order equations.} In this section we apply the $\sd$-Galois theory to understand the differential properties of solutions of higher order linear difference equations. The main tool will be the following proposition.

\begin{prop}\label{higher_order_prop} Let $k$ be a $\sd$-field with $k^\sigma$ $\der$-differentially closed and $A \in \GL_n(k)$.  Assume that the $\sigma$-Galois group of $\sigma(Y) = AY$ is a simple, noncommutative linear algebraic group of dimension $t$. Let $R = k\{Z,\frac{1}{\det Z}\}_\der$ be the $\sd$-PV ring for this equation.  The differential dimension of $R$ over $k$ is less than $t$ if and only if there exists a $B \in \gl_n(k)$such that $\sigma(B) = ABA^{-1} + \der(A)A^{-1}$ (in which case $Z^{-1}(\der(Z) - BZ) \in \gl_n(k^\sigma)$ and the differential dimension of $R$ is $0$).
\end{prop}
\begin{proof} Let $B \in \gl_n(k)$ satisfy $\sigma(B) = ABA^{-1} + \der(A)A^{-1}$.  A computation shows that $\sigma(\der(Z) - BZ) = A(\der(Z) - BZ)$.  Therefore, the columns of $\der(Z) - BZ$ are solutions of $\sigma(Y) = AY$ implying that there exists a $D \in \gl_n(k^\sigma) $ such that $\der(Z) - BZ = ZD$ or $Z^{-1}(\der(Z) - BZ) \in \gl_n(k^\sigma)$. Since $\der(Z) = BZ + ZD$, we have that each entry of $Z$ is differentially algebraic over $k$ and so the differential dimension of $R$ over $k$ is zero.

Now assume that the differential dimension of $R$ over $k$ is less than $t$.  Theorem~\ref{torsor} implies that the differential dimension of $R$ is equal to the differential dimension of the $\sd$-Galois group $H$ (see Proposition~\ref{dimprop}).  The differential dimension of $G$, considered now as a differential group, is again $t$ (\cite{DAAG}, Proposition 10, p.~200) so we can conclude from Proposition~\ref{sdzariskidense} that $H$ is a Zariski-dense {\em proper} differential subgroup  of $G(k^\sigma)$.  In \cite{cassidy6}, Phyllis Cassidy showed that such a group is conjugate over $k^\sigma$ to $G(C)$, where $C = \{c \in k^\sigma \ | \der(c) = 0\}$.  Proposition~\ref{constantgaloisgroup} implies the conclusion of this proposition. \end{proof}

Once again the assumption that $k^\sigma$ is differentially closed is rather artificial but this result can nonetheless be used to prove results about meromorphic functions. The following corollary follows from Proposition~\ref{higher_order_prop} in much the same was as Corollary~\ref{indefsumcor} follows from Proposition~\ref{indefsum} and its proof will be left to the reader.

\begin{cor} \label{higher_order_cor} Let $k$ be a $\sd$-field satisfying hypothesis (A) or (B) above.  Let $A(x) \in \GL_n(\CX(x))$ and assume
 that the difference equation $\sigma(Y) = AY$ has $\sigma$-Galois group $\SL_n$ over $\CX(x)$.   Let $Z$ be a fundamental solution matrix of this equation, meromorphic in $\CX$ if (A) holds or in $\CX^*$ if (B) holds.  Then the differential dimension of $k\{Z, \frac{1}{\det Z}\}_\der$ over $k$ is less
 than $n^2-1$ if and only if there exists a $B \in \gl_n(\CX(x))$ such that $\sigma(B) = ABA^{-1} + \der(A)A^{-1}$.
\end{cor}
Note that the system $\sigma(B) = ABA^{-1} + \der(A)A^{-1}$ is an inhomogeneous system of linear difference equations in the entries of $B$ and that there are well developed algorithms to determine if such a system has a solution whose entries are rational functions, \cite{abramov_zima, barkatou99, bron_diff, hoeij98b, hoeij_sing, PWZ,  wu05}.  We will give two examples to illustrate {\em ad hoc} strategies.

\begin{ex}{\em In \cite{PuSi}, p.~42, it is shown that the difference equation $Y(x+1) = A(x)Y(x)$, where \[A(x) = \left(\begin{array}{cc} 0 &-1\\1&x \end{array}\right) \]
has Galois group over $\CX(x)$ equal to $\SL_2(\CX)$.  We shall show that the corresponding differential dimension of the $\sd$-PV extension is $3$.  Assuming this is not the case,  Corollary~\ref{higher_order_cor} implies that there exists a $B \in \gl_n(\CX(x))$ such that $\sigma(B) = ABA^{-1} + \der(A)A^{-1}$. Letting $B = \left(\begin{array}{cc} a &b\\c&d\end{array}\right)$, this implies that
\[\sigma\left(\begin{array}{c} a\\b\\c\\d \end{array}\right) =
\left(\begin{array}{cccc} 0&0&-x&1\\0&0&-1&0\\x&-1&x^2&-x\\1&0&x&0 \end{array}\right)\left(\begin{array}{c} a\\b\\c\\d \end{array}\right) + \left(\begin{array}{c} 0\\0\\-1\\0 \end{array}\right)\]
where $\sigma(x) = x+1$. Eliminating the $a,c$ and $d$ one sees that $b$ satisfies
\begin{eqnarray}\label{shiftex}x\sigma^3(b)  -(x^3+2x^2-1)\sigma^2(b) + x(x^2+x-1)\sigma(b) -(x+1)b=
2x+1\end{eqnarray} (this can be done by hand or using the Maple
commands given in \cite{HaSi}). We now claim that this latter
equation has no rational solution.  We first show that any rational
solution must be a polynomial.  If $\alpha$ is a pole of a putative
solution $b(x)$, let $n$ be the largest nonnegative integer such
that $\beta = \alpha-n$ is again  a pole of $b(x)$.  If $\beta$ is
not an integer then $\sigma^3(b)$ has a pole at $\beta -3$ that is
not cancelled in $(\ref{shiftex})$.  Therefore all possible poles of
$b$ must be at integers.  If $\lambda$ is  the largest integer pole
then the poles of $\sigma(b), \sigma^2(b)$ and $\sigma^3(b)$ are all
smaller than $\lambda$. Therefore, to be cancelled in
$(\ref{shiftex})$, $\lambda = -1$.  If $\mu$ is the smallest pole,
then $\mu-3$ is a pole of $\sigma^3(b)$ that is not cancelled unless
$\mu-3 = 0$.  Since $\mu \leq \lambda$ we have a contradiction
unless $b$ is a polynomial.  Now assume that $b=
b_0x^n+b_1x^{n-1}+b_2x^{n-2} + b_3x^{n-3}+\ldots$ is a polynomial
and substitute in the above equation. On the left hand side of the
equal sign  the coefficient of $x^{n+3}$ is zero and the coefficient
of $x^{n+2}$ is $-(n+1)b_0$.  Since $n\geq 0$, this term must be
zero (or it would not cancel with a term on the right hand side of
the equation)  and we get a contradiction.}
\end{ex}
\begin{ex}\label{q-diffex}{\em Difference Galois groups of $q$-hypergeometric equations have been calculated by Hendriks \cite{hendriks_qdiff} and Roques \cite{roques}.  We will consider the following family of equations from Roques classification:
\[y(q^2x) - \frac{2ax -2}{a^2x-1}y(qx) + \frac{x-1}{a^2x-q^2}y(x) = 0
\]
where $a \notin q^{\ZX}$ and $a^2 \in q^\ZX$.  Roques showed (see Section 4.2, Theorem 10 of \cite{roques}) that equations of this form have difference Galois group $\SL_2(C)$.  In \cite{HaSi}, we apply Corollary~\ref{higher_order_cor} to show that the associated $\sd$-PV extension has differential transcendence degree 3, that is that, for $y_1(x), y_2(x)$  solutions linearly independent over the field $\calG$ of $q$-periodic functions,  $y_1(x), y_2(x), y_1(qx)$ are differentially independent over $\calG(x)$.  The computations of \cite{HaSi} give rise to large formulas and so we will not present the proof of full result here.  In this paper we will present the arguments for the case that $q = 1/4$ and $a=1/2$ and note that the arguments given in \cite{HaSi} are just simple generalizations of the arguments given here.\\[0.2in]
When $q = 1/4$ and $a=1/2$, we are considering the equation
\[y(q^2x) - \frac{4(x -2)}{x-4}y(qx) + \frac{16(x-1)}{4x-1}y(x) = 0.
\]
We wish to show that the equation $B(qx) = ABA^{-1} +
x\frac{d}{dx}(A)A^{-1}$ has no solution $B \in \gl_2(\CX(x))$.
Letting $B = \left(\begin{array}{cc} u &v\\w&z\end{array}\right)$,
we show in \cite{HaSi} that $v$ satisfies the third order
equation{\small
\begin{eqnarray*}\label{qhypex}v(q^3x)
- 1/4\,{\frac { \left( x-64 \right)  \left( x-4 \right)  \left( 20\,{x}^
{2}-353\,x+1032 \right) }{ \left( x-32 \right)  \left( -1+x \right) ^{
2} \left( x-16 \right) }}v(q^2x)
+ {\frac { \left( 20\,{x}^{2}-353\,x+1032 \right)  \left( x-64 \right)
}{ \left( x-16 \right)  \left( 4\,x-1 \right)  \left( x-32 \right) }}v(qx)\\
+1/4\,{\frac { \left( x-64 \right)  \left( x-2 \right)  \left( 4\,x-1
 \right) }{ \left( -1+x \right) ^{2} \left( x-32 \right) }}v(x) =
-{\frac {x \left( x-64 \right)  \left( 47\,{x}^{2}-496\,x+1952
 \right) }{ \left( 4\,x-1 \right)  \left( x-8 \right)  \left( x-16
 \right)  \left( -1+x \right)  \left( x-32 \right)}}
\end{eqnarray*}}

\noindent The key observation is that $8 = 1/(qa)$ is a pole of the right hand side of this equation but is not a pole of any of the coefficients on the left hand side. Therefore $8$ must be a pole of some $y(q^ix), i=0,1,2,3$.  Before we proceed, we make two definitions.  If $R(x)$ is a rational function, we say that $\alpha$ is a {\em maximal pole} of $R(x)$ if $\alpha$ is a pole of $R$ but $\alpha/q^n$  is not a pole of $R$ for any $n>0$ and  $\beta$ is a {\em minimal pole} of $R(x)$  if $\beta$ is a pole but $\beta q^n$  is not a pole of $R$ for any $n>0$.  We shall now show that the assumption that $x=8$ must be a pole of some $y(q^ix), i=0,1,2,3$ leads to a contradiction.  If $x=8$ is a pole of $v(x)$, then $v(x)$ has a maximal  pole of the form $8(4^n)$ for some $n \geq 0$.  In this case $v(q^3x)$ will have a maximal pole of the form $8(4^{n+3}) \geq 512$ and this could  not possibly cancel with anything in the above equation.   A similar argument shows that $v(qx)$ cannot have a pole at $8$. If $v(q^2x)$ has a pole at $x=8$, then it has a minimal pole of the form $8/4^n$.  Therefore $v(x)$ will have a minimal pole of the form $8/4^{n+2}$. No other term in the equation has a pole of this form so there again can be no cancelation.  A similar argument shows that $v(q^3x)$ has no pole at $x= 8$.  Therefore the above equation cannot have a rational solution.
}
\end{ex}

\section{Inverse Problem}\label{inverse}

In this section we consider $\sd$-fields of the form $k=C(x)$ where $C$ is a differentially closed $\der$-field containing $\CX$, $x$ is transcendental over $C$ and either $\sigma(x) = x+1$, $\der(x) = 1, \ \sigma|_C = id$ or $\sigma(x) = qx, \ |q|\neq 1$, $\der(x) = x, \ \sigma|_C = id$.  We will consider the following special cases of the inverse problem:

\begin{quotation} \noindent {\em Which differential subgroups of the additive group $\GX_a(C)$ are $\sd$-Galois groups of $\sd$-PV extensions of $k$ and which of these are defined by equations of the form $\sigma(y) - y = f, \ f\in \CX(x)$}?
\end{quotation}

The answer to the first question (Proposition~\ref{generalinverseprop}) is that all such subgroups can occur and to the second (Proposition~\ref{inverseprop}) is that very few occur.  We begin by showing that a $\sd$-Picard-Vessiot ring with $\sd$-Galois group a subgroup of $\GX_a$ is generated by a solution of an equation of the form $\sigma(y)-y = f$.

\begin{lem} Let $k = C(x)$ be as above and let $R$ be a $\sd$-PV extension with $\sd$-Galois group $G$ isomorphic to  a differential subgroup of $\GX_a(C)$.  Then $R = k\{t\}_\der$ where $\sigma(t) - t \in k$.

\end{lem}
\begin{proof} We may write $R = k\{Z, \frac{1}{\det Z}\}_\der$ where $\sigma(Z)Z^{-1} \in \GL_n(k)$ and $G \subset \GL_n(C)$.    From Proposition 35 of \cite{cassidy1}, we know that $G$ is unipotent. From Proposition~\ref{sdzariskidense} we have that the Zariski closure $H$ of $G$ is the $\sigma$-Galois group of the $\sigma$-PV extension $S = k[Z, \frac{1}{\det Z}]$.  Since $H$ is unipotent and commutative, we have that $H \simeq (\GX_a)^m$ as an algebraic group.  This latter group has first Galois cohomology set that is trivial and so the $H$-torsor corresponding to $S$ is also trivial, implying that $S = k[t_1, \ldots , t_m]$ where  for any $\phi \in H(C)$, $\phi(t_i) = t_i +c_{i,\phi}$ for some $c_{i,\phi} \in C$. Note that this latter fact implies that $\sigma(t_i) - t_i \in k$.

Now assume that $G \simeq \GX_a(C)$ as a linear differential algebraic group.  Each of the maps $\Phi_i: \phi \mapsto c_{i,\phi}$ defines a differential homomorphism from $\GX_a(C)$ to $\GX_a(C)$ and therefore, identifying $\phi \in G$ with an element $d_\phi \in \GX_a(C)$, $\Phi_i$ must be of the form $\Phi_i(d_\phi) = L_i(d_\phi)$ for some $L_i \in C[\der]$ (Corollary 4 of \cite{cassidy2}).  The map from $\GX_a(C)$ to $\GX_a(C)$ given by  $d\mapsto (L_1(d), \ldots , L_m(d))$ is an isomorphism and so has trivial kernel.  Since $C$ is differentially closed, we must have that the equations $L_1(y) = \ldots = L_m(y) = 0$ have no nonzero solution and so the left ideal $C[\der]L_1+\ldots C[\der]L_m \subset C[\der]$ contains $1$. Let $N_1L_1+\ldots +N_mL_m = 1, \ N_i \in C[\der]$ and let $t = N_1(t_1) + \ldots +N_m(t_m) \in R$.  For any $\phi \in G$, we have $\phi(t) = \phi(\sum_i N_i(t_i)) = \sum_iN_i(t_i + L_i(d_\phi)) = \sum_i N_i(t_i) + (\sum_iN_iL_i)(d_\phi) = t+d_\phi$.  Therefore, if $\phi \neq id$, then $\phi(t) \neq t$.  This implies that the only automorphism that leaves $k\{t\}_\der$ fixed is the identity.  From the Galois correspondence we know that $R = k\{t\}_\der$.

Now assume that $G$ is isomorphic to a proper differential subgroup of $\GX_a$.  In this case the differential dimension of $G$ is zero and furthermore $G$ is connected.  Therefore $R$ is domain.  The Differential Primitive Element Theorem \cite{seidenberg52} implies that there exist $c_i \in C$, not all zero such that $R = k\{t\}_\der$ where $t = \sum_ic_it_i$.  Clearly, $\sigma(t) - t = \sum_ic_i(\sigma(t_i) - t_i) \in k$.\end{proof}

\noindent The above result reduces the problem of finding which differential subgroups of $\GX_a$ occur as $\sd$-Galois groups over $k$ to the problem of deciding which subgroups of $\GX_a$ occur as $\sd$-Galois groups of equations of the form $\sigma(y) - y = f, \ f\in k$.

\begin{prop}\label{generalinverseprop} Let $k=C(x)$ be a $\sd$-field as above and $G \subset \GX_a(C)$ a linear differential algebraic group.  Then there exists an $f \in C(x)$ such that $\sigma(y) - y = f$ has $\sd$-PV group $G$ over $k$.
\end{prop}
\begin{proof}  Let us assume that $G$ is a proper subgroup of $\GX_a$.  In this case, there is a nonzero  linear operator $L \in C[\der]$ such that $G = \{c \in C \ | \ L(c) = 0 \}$.  Since $C$ is differentially closed, there exist $c_1, \ldots , c_n \in C$, linearly independent over $C^\der$, whose $C^\der$-span is the solution space of $L(y) = 0$.

Now assume that $\sigma(x) = x+1$ and $\der(x) = 1$. We can also deduce from $C$ being differentially closed, that there is an $x_0 \in C$ such that $\der(x_0) = 1$ and therefore that $\der(x-x_0) = 0$. Consider the equation
\[\sigma(y) - y = \sum_{i=1}^n \frac{c_i}{(x-x_0)^i} \ .\]
Note that  $\sigma(L(y)) - L(y) = L(\sigma(y) - y) = L(\sum_{i=1}^n \frac{c_i}{(x-x_0)^i}) = \sum_{i=1}^n \frac{L(c_i)}{(x-x_0)^i} = 0$.  Therefore $L(y) \in C$.
We claim that the $\sd$-Galois group $H$ of this equation is $G$.  To see that $H \subset G$, let $\phi \in H$. There exists a $d \in C$ such that $\phi(y) = y+d$.  Since $L(y) \in C$, we have that $L(d) = \phi(L(y)) - L(y) = 0$.  Therefore $H \subset G$.

To see that $G\subset H$, let $H = \{c \ | \ \tilde{L}(c) = 0\}$.  A calculation shows that $\tilde{L}(y)$ is left fixed by $H$ and so $\tilde{L}(y) = g \in k$.  This implies that \[\tilde{L}(\sigma(y) - y) = \tilde{L}(\sum_{i=1}^n \frac{c_i}{(x-x_0)^i}) = \sum_{i=1}^n \frac{\tilde{L}(c_i)}{(x-x_0)^i}= \sigma(g) - g\ . \]
Since the polar dispersion of $\sum_{i=1}^n \frac{\tilde{L}(c_i)}{(x-x_0)^i}$ is zero, we must have, by Lemma~\ref{posdisp} and the uniqueness of partial fraction decompositions, that each $\tilde{L}(c_i) = 0$.  Therefore $G \subset H$.

Now assume that $\sigma(x) = qx$ and $\der(x) = x$.  We can also deduce from $C$ being differentially closed, that there is an $x_1 \in C$ such that $\der(x_1) = x_1$ and therefore that $\der(\frac{x}{x_1}-1) = 0$. Consider the equation
\[\sigma(y) - y = \sum_{i=1}^n \frac{c_i}{(\frac{x}{x_1}-1)^i}\]
A similar argument as above shows that this equation has $\sd$-Galois group $G$.

We have shown that every proper subgroup of $\GX_a$ is a $\sd$-Galois group.  Examples \ref{ex1} and \ref{ex2} show that $\GX_a$ can be realized as well. \end{proof}

We now turn to the problem of characterizing which subgroups of $\GX_a(C)$ occur as $\sd$-Galois groups of equations of the form $\sigma(y) - y = f$ with $f \in \CX(x)$.

\begin{prop}\label{inverseprop}Let $C$ be a differentially closed $\der$-field containing $\CX$ and $k=C(x)$, $x$ transcendental over $C$. \begin{enumerate}
\item Assume $\sigma(x) = x+1, \der(x) = 1, \sigma|_C = id, \der |_C = \der$.  A differential subgroup $G$ of $\GX_a$ is the $\sd$-Galois group of an equation of the form $\sigma(y) - y = f, \ f \in \CX(x)$ if and only if
\begin{enumerate}
\item $G = \{0\}$, when $f = \sigma(g) - g$ for some $g \in \CX(x)$,or
\item $G = \GX_a(C)$, in all other cases.
\end{enumerate}
\item Assume $\sigma(x) = qx, \ |q| \neq 1, \der(x) = x, \sigma|_C = id, \der |_C = \der$. A differential subgroup $G$ of $\GX_a$ is the $\sd$-Galois group of an equation of the form $\sigma(y) - y = f, \ f \in \CX(x)$ if and only if
\begin{enumerate} \item $G = \{0\}$, when $f = \sigma(h) - h$ for some $h\in \CX(x)$, or
\item  $G = \GX_a(C^\der)$, when $f = \sigma(h) - h + c $ for some $h\in \CX(x)$ and nonzero $c \in \CX$, or
\item $G = \GX_a(C)$ in all other cases.
\end{enumerate}
\end{enumerate}
\end{prop}
\begin{proof} Let $R = k\{z\}_\der$ be a $\sd$-PV extension of $k$ with $\sigma(z) - z = f$. Assume that the $\sd$-Galois group is a proper subgroup of $\GX_a$ and let $G = \{c \ | \ L(c) = 0\}$ for some nonzero $L \in C[\der]$.  One sees that $L(z)$ is left invariant by $G$ so $L(z) = g \in k$.  We then have that $L(f) = \sigma(L(z)) - L(z) = \sigma(g) - g$.  The coefficients of powers of $\der$ in $L$ and the coefficients of the powers of $x$ in $g$ are in $C$ while the coefficients of $f$ are in $\CX$.  Equating powers of $x$ in $L(f) = \sigma(g) - g$, we get a system of polynomial equations, over $\CX$, for the coefficients of powers of $\der$ in $L$ and the coefficients of the powers of $x$ in $g$.  Since these have a solution, we can find a solution in $\CX$.  Therefore we may assume that we have a relation of the form $L(f) = \sigma(g) - g$ where $L$ is a nonzero element of $\CX[\der]$ and $g \in \CX(x)$.

In case 1.~above,  Lemma~\ref{trivlem} implies that $f = \sigma(h) - h$ for some $h \in \CX(x)$.  Therefore, $\sigma(y) - y = f$ has a solution in $\CX(x)$ and so has trivial Galois group.

In case 2.~above, Lemma~\ref{trivlem} implies that $f = \sigma(h) - h +c$ for some $h \in \CX(x), c\in \CX$. If $c = 0$ then the Galois group is trivial.  If $c \neq 0$, then the $\sd$-Galois group $G$ cannot be trivial. Furthermore, we have that $\der(f)$ satisfies $\sigma(\der(f)) - \der(f) = \sigma(\der(h)) -\der(h)$ and so $\der(f) \in k$.  Let $\phi \in H$ and assume that $\phi(f) = f+d$.  We then have that  $\der(f)= \phi(\der(f)) = \der(f) + \der(c)$ so $\der(c) = 0$.  Therefore $G \subset \CX_a(\CX)$.  Since the only proper differential subgroup of this group  is $\{0\}$ and we have excluded this, we must have $G = \GX_a(\CX)$. Examples of each case are given below.\end{proof}
\begin{ex}\label{ex1} Let $k = C(x)$ be the $\sd$-field where $\sigma(x) = x+1, \der(x) = 1, \sigma|_C = id, \der |_C = \der$.  If $f = \frac{1}{x}$, the $f$ has a pole and its polar dispersion is $0$ so Lemma~\ref{posdisp} implies that $f \neq \sigma(g) - g$ for any $g \in \CX(x)$.  Therefore,
\[y(x+1) - y(x) = \frac{1}{x}\]
has $\sd$-Galois group $\GX_a(C)$.
\end{ex}
\begin{ex}\label{ex2} Let $k = C(x)$ be the $\sd$-field where $\sigma(x) = qx, \ |q| \neq 1, \der(x) = x, \sigma|_C = id, \der |_C = \der$. If $f = \frac{1}{x-1}$, then an argument similar to the one in Example~\ref{ex1}, shows that
\[y(qx) - y(x) = \frac{1}{x-1}\]
has $\sd$-Galois group $\GX_a(C)$.  The above proposition implies that $y(qx) - y(x) = qx-x+1$ has $\sd$-Galois group $\GX_a(\CX)$.
\end{ex}

\section{Parameterized Difference Equations}\label{param}  In the previous sections, we applied the $\sd$-PV theory to the field $\CX(x)$ with the derivation and automorphism acting nontrivially on the same variable, $x$. In this section we will consider  the $\sd$-PV theory applied to study  parameterized difference equations $\sigma(Y) = A(x,t) Y$ where $\sigma$ acts only on the variable $x$ and the derivation we consider is the partial derivative with respect to $t$. We shall only consider one small aspect of this: how the $\sd$-Galois group measures to what extent  the connection matrix of a parameterized system of regular $q$-difference equations actually depends on $t$.  It is hoped that this will stimulate further research along these lines. \\[0.2in]
 We shall consider parameterized difference equations.  Let $\calO$ be an open subset of the complex plane and let $\calM$ be the field of functions $f(t)$  meromorphic on $\calO$. We shall consider difference equations of the form
\begin{eqnarray}\label{parameq}
Y(qx,t) = A(x,t)Y(x,t)
\end{eqnarray}
where $A(x,t) \in \GL_n(\calM(x))$ and $A(0,t)= I_n = A(\infty, t)$, independent of $t$ and $|q| \leq 1$. These are a parameterized version of the equations considered by Etingof in  \cite{etingof}. Following Etingof, one sees that
\[Y_0(x,t) = \prod_{j=0}^\infty A(q^jx,t)^{-1}  = A(x,t)^{-1}A(qx,t)^{-1} \cdots \]
is  a solution of (\ref{parameq}), meromorphic on $\CX\times \calO$   and
\[Y_{\infty}(x,t) = \prod_{j=1}^\infty A(q^{-j}x,t) = A(x,t)A(q^{-1}x,t) \cdots \]
is  a solution of (\ref{parameq}), meromorphic on $(\CX^*\cup\infty) \times \calO$. \\[.2in]
Let $k_0$ be a field containing $\calM$, differentially closed with
respect to $\dd = \frac{d}{dt}$.  We note that any differential
subfield $  k_1 \subset k_0$ that is finitely generated over $\QX$
in the differential sense, may be identified with a differential
field of functions  meromorphic on some subdomain of $\calO$
\cite{sei58,sei69}. In particular we may consider any finite set of
elements of $k_0$ as being functions on some domain.   Let $k =
k_0(x)$. Finally we let $C(x,t) = Y_0(x,t)^{-1}Y_\infty(x,t)$. Note
that the entries of $C(x,t)$ are meromorphic in $t$ and $q$-periodic
meromorphic  in $x$. In \cite{etingof} Etingof proves the
non-parameterized versions of the following facts and his proofs
readily generalize to the present case:
\begin{enumerate}
\item Each of the fields $L_0 = k<Y_0>_\dd$ and $L_\infty = k<Y_{\infty}>_\dd$ are $\sd$-PV extensions of $k$.
\item Let $w\in \CX^*$ be such that $C(w,t)$ is defined and  invertible. Then the map $Y_\infty \mapsto Y_0C(w,t)$ defines an isomorphism of $\sd$-PV extensions $\tau_w:L_\infty\rightarrow L_0$. Therefore, if $u, w$ are values where $C$ is defined and invertible, then $C(u,t)C(w,t)^{-1} \in \Aut_\sd (L_\infty/k)$
\item If $\Gamma$ is the subgroup of $\GL_n(k_0)$ generated by all matrices $C(u,t)C(w,t)^{-1}$, where $u, w$ are values where $C$ is defined and invertible, then the fixed field of $\Gamma$ is $k$.
\end{enumerate}
We show how Etingof's ideas can be modified to prove statement 3.~above.  Let $h(x,t) \in L_\infty$ be fixed by all matrices $C(u,t)C(w,t)^{-1}$.  We have that $h = P(x,t,Y_\infty(x,t)C(u,t)C(w,t)^{-1})$, where $P$ is a rational expression in $x$ and $Y_\infty(x,t)C(u,t)C(w,t)^{-1}$ and its derivatives with respect to $t$.  Since this holds for an uncountable number of $u \in \CX$, we may replace $u$ by $x$ and we have $h(x,t) = P(x,t,Y_\infty(x,t)C(x,t)C(w,t)^{-1}) = P(x,t,Y_0(x,t)C(w,t)^{-1})$.  This shows that, for all values of $t$ in some subdomain of $\calO$ we have that $h(x,t)$ is  meromorphic on the Riemann Sphere.  This implies that for each such value of $t$, $h$ is a rational function of $x$. We wish to show that $h(x,t)$  it is a rational function of $x$ with coefficients in $k_0$.  Note that for some uncountable set $S$ of values of $t$, there exists an integer $N$ such that for $t_1 \in S$, $h(x,t_1)$ is the quotient of polynomials of degree $N$ in $x$.  This implies that the Casoratian $Cas(h(x,t), xh(x,t), \ldots , x^Nh(x,t), 1, x, \ldots , x^N)$ vanishes for an uncountable  set of values of $t$ and so must be identically zero.  Therefore $h(x,t), xh(x,t), \ldots , x^Nh(x,t), 1, x, \ldots , x^N$ are linearly dependent over $k_0$ and we have that $h(x,t)$  it is a rational function of $x$ with coefficients in $k_0$.  \\[0.2in]
The Galois correspondence therefore allows us to conclude the following result.
\begin{prop} For $L_\infty$and $\Gamma$ as above, we have that the Kolchin closure of $\Gamma$ is $\Aut_\sd (L_\infty/k)$. \end{prop}
In \cite{etingof}, Etingof shows how his similar result can be used to deduce that any connected linear algebraic group is the Picard-Vessiot group of a $q$-difference equation over $\CX(x)$.  It would be interesting to see to what extent the above result can be used to attack the inverse problem over $\CX(x,t)$. \\[0.2in]
In analogy to the nonparametric case, we shall refer to $C(x,t)$ as
the  {\em parameterized connection matrix}.  The following
proposition states that the parametric connection matrix is
equivalent  to a  matrix independent of $t$ if and only if the
$\sd$-Galois group is conjugate to a constant group.
\begin{prop}  Let $k, A(x,t), C(x,t)$ be as above.  The following are equivalent:
\begin{enumerate} \item The $\sd$-Galois group $\Aut_\sd(L_\infty/k)$  is conjugate over $k_0$ to a subgroup of $\GL_n(\CX)$.
\item There exists  matrices $D, E\in \GL_n(k_0)$ such that $\frac{d}{dt}(DCE) = 0$.
\item There exist  solutions $\bar{Y}_0(x,t)$, analytic at $x=0$, and $\bar{Y}_\infty(x,t)$, analytic at $x=\infty$,    of $Y(qx,t) = A(x,t)Y(x,t)$ such that $\bar{C} = \bar{Y}_0^{-1}\bar{Y}_\infty$ satisfies $\frac{d}{dt}(\bar{C}) = 0$.
\item There exists a matrix $B \in \gl_n(k)$ such that $\sigma(B) = AB A^{-1} + \dd(A)A^{-1}$, in which case the system
\begin{eqnarray*}
\sigma(Y)  & = & AY\\
\dd(Y) & = & BY
\end{eqnarray*}
has solution $Y = Y_\infty D$ for some $D \in\GL_n(k_0)$.
\end{enumerate}
\end{prop}
\begin{proof} We shall first show that $1$.~is equivalent to $2$. Let $D(t) \in \GL_n(k_0)$ conjugate $\Aut_\sd(L_\infty/k)$ into $\GL_n(\CX)$.  We then have that $D(t)C(u,t)C(w,t)D^{-1}(t) \in \GL_n(\CX)$ for an uncountable number of values of $u\neq w \in \CX$.  Let $H(x,t) = D(t)C(x,t)$.  Differentiating $H(u,t)H^{-1}(w,t)$ with respect to $t$, one sees that \[H^{-1}(u,t)\frac{d}{dt}H(u,t) = H^{-1}(w,t)\frac{d}{dt}H(w,t)\] for an uncountable number of $u\neq w$.  This implies that $H^{-1}(x,t)\frac{d}{dt}H(x,t) = J(t) \in \gl_n(k_0)$.  Since $k_0$ is differentially closed, there exists a $Z(t) \in \GL_n(k_0)$ such that  $\frac{d}{dt}Z(t) = Z(t) J(t)$ and so $\frac{d}{dt}(H(x,t)Z^{-1}(t)) = 0$.  Setting $E(t) = Z^{-1}(t)$ gives us the conclusion of $2$.\\[0.2in]
Assuming $2.$, one sees that $\frac{d}{dt}(D(t)C(u,t)E(t)E^{-1}(t)C^{-1}(w,t)D^{-1}(t)) = 0$. Since the elements $C(u,t)E(t)E^{-1}(t)C^{-1}(w,t)$ generate a Kolchin dense subgroup of $\Aut_\sd(L_\infty/k)$, we have that $D$ conjugates this group into a subgroup of $\GL_n(k_0)$. Therefore $1$.~holds.\\[0.2in]
To see that $2$.~is equivalent to $3$., assume that there exists a matrices $D, E\in \GL_n(k_0)$ such that $\frac{d}{dt}(DCE) = 0$.  Letting $\bar{Y}_0 = Y_0 D^{-1}$ and $\bar{Y}_\infty = Y_\infty E$ yields elements that satisfy $3$.  Conversely, given $\bar{Y}_0$ and $\bar{Y}_\infty$ as above, we have that there exist $D^{-1}, E\in \GL_n(k_0)$ such that $\bar{Y}_0 = Y_0D^{-1}$ and $\bar{Y}_\infty = Y_\infty E$.  This implies that $\frac{d}{dt}(DCE) = 0$.\\[0.2in]
The equivalence of $1$.~and $4$.~follows from Proposition~\ref{constantgaloisgroup}.
\end{proof}
We note that the reasoning of this last result can be adapted to the more general situation of a parameterized system of $q$-difference equations that are regular at $0$ and $\infty$ ({\em cf.,} \cite{PuSi}, Ch. 12.3.2) as well as shift difference equations that are regular at $\infty$ ({\em cf.,} \cite{PuSi}, Ch. 8.5).  It would be interesting to develop analogous results for parameterized equations with  singularities.

\section{Appendix}
\subsection{Rational Solution of Difference Equations}\label{diffappendix}  In this section we have collected some elementary facts concerning normal forms and solutions of difference equations that are needed in the previous sections. Many of these results appear implicitly (and some explicitly) in the literature \cite{abramov75, bronstein92,karr, matusevich, paule2, PuSi,schneider_phd, schneider_bounds, schneider_telescope} but we assemble them here in a form that meets our needs.    We begin by considering the field $C(x)$, $C$ algebraically closed, with the difference operator $\sigma(x) = x+1$ or $\sigma(x) = qx$, $q$ not a root of unity.

\begin{defin}\label{dispdef}({\em cf.}, \cite{abramov75, paule2}) Let $k = C(x)$ be as above.  $P,Q \in C[x], f = \frac{P}{Q}, \ \gcd(P,Q) = 1$. \begin{enumerate}
\item If $\sigma(x) = x+1$, we define the {\em dispersion of $Q$, $\disp(Q)$} to be the largest nonnegative integer $\ell$ such that for some $\alpha \in C$, $\alpha$ and $\alpha + \ell$ are roots of $Q$.

\item If $\sigma(x) = qx$, we define the {\em dispersion of $Q$, $\disp(Q)$} to be the largest nonnegative integer $\ell$ such that for some {\em nonzero} $\alpha \in C$, $\alpha$ and $q^\ell\alpha$ are roots of $Q$.
\item We define the {\em polar dispersion of $f$, $\pdisp(f)$} to be the dispersion of $Q$.
\item The element $f$ is said to be {\em standard} if $\disp(P\cdot Q) = 0$ ({cf.}, \cite{PuSi} Ch.2).

\end{enumerate}

\end{defin}

\begin{lem}\label{normalform} Let $k= C(x)$ be as above, $f \in k$ and $a \in C, a \neq 0$.  \begin{enumerate}
\item There exist $f^*,g \in k$ with  $\pdisp(f^*) = 0$   such that $f = f^* + \sigma(g) - ag$.
\item Assume $f \in k\backslash\{0\}$. There exist $\tilde{f}, \tilde{g}\in k\backslash\{0\}$ such that $f = \tilde{f}\frac{\sigma{\tilde{g}}}{\tilde{g}}$ with $\tilde{f}$ standard ({\em cf.}~Lemma 2.2, \cite{PuSi}).
\end{enumerate}
\end{lem}
\begin{proof} 1.~Let us first assume that $\sigma(x) = x+1$.
For any $f = \frac{P}{Q}, \ \gcd(P,Q) = 1$, let $n_f $ be the number
of distinct roots $\alpha$ of $Q$ such that $\alpha +N$ is also a
root of $Q$, where $N = \disp(Q)$.  We order the pairs $(\disp(f),
n_f)$ lexicographically and proceed to prove the lemma by induction.
If $\pdisp(f) = 0$ we are done.  Let $\pdisp(f) = N>0$ and $n_f >0$.
Let $\alpha \in C$ be a root of $Q$ such that $\alpha - N$ is again
a root of $Q$. We may write
\[ f = \sum_i \frac{a_i}{(x-(\alpha-N))^i} + h\]
where $\alpha-N$ is not a root of the denominator of $h$. Let $\tilde{g} = \sum \frac{a_i}{x-(\alpha-(N-1))^i}$ and let $\tilde{f} = f- (\sigma(\tilde{g}) + a\tilde{g})$.  One sees that either $\pdisp(\tilde{f}) < \pdisp(f)$ or $\pdisp(\tilde{f}) = \pdisp(f)$ and $n_{\tilde{f}} < n_f$.  Therefore $\tilde {f} = f^* +\sigma(\bar{g})-\bar{g}, \ \pdisp(f^*) = 0$. This implies that $f = f^* + \sigma(\tilde{g} + \bar{g}) - a(\tilde{g} + \bar{g})$.\\[0.1in]
Now assume we are in the case where $\sigma(x) = qx$.  For any $f =
\frac{P}{Q}, \ \gcd(P,Q) = 1$, let $n_f $ be the number of distinct
{\em nonzero} roots $\alpha$ of $Q$ such that $q^{-N}\alpha $ is
also a root of $Q$, where $N = \disp(Q)$.  We order the pairs
$(\disp(f), n_f)$ lexicographically and proceed to prove the lemma
by induction. If $\pdisp(f) = 0$ we are done.  Let $\pdisp(f) = N>0$
and $n_f >0$. We may assume $C$ is algebraically closed. Let $\alpha
\in C$ be a root of $Q$ such that $q^{-N}\alpha$ is again a root of
$Q$. We may write
\[ f = \sum_i \frac{a_i}{(x-q^{-N}\alpha)^i} + h\]
where $q^{-N}\alpha$ is not a root of the denominator of $h$. Let $\tilde{g} = \sum_i \frac{q^ia_i}{(x-q^{-(N-1)}\alpha)^i}$.  We then have that $\sigma(\tilde{g}) = \sum_i \frac{q^ia_i}{(qx-q^{-(N-1)}\alpha)^i} = \sum_i \frac{a_i}{(x-q^{-N}\alpha)^i}$ so we can apply the induction assumption to
$\tilde{f} = f- a(\sigma(\tilde{g}) - \tilde{g})$, and achieve the conclusion of the lemma as before. \\[0.1in]
2.~ Assume that $\sigma(x) = x+1$.  Let $\dd = \frac{d}{dx}$ and
\[ h = \frac{\dd(f)}{f} = \sum_i\frac{n_i}{x-\alpha_i}\] for some integers $n_i$ and $\alpha_i \in C$.  If one applies the  method  of part 1.~to $h$, one finds that $h = h^* + \sigma(g) - g$, where $\pdisp(h^*) = 0$ and $h^* = \sum \frac{m_i}{x-\beta_i}, \ g = \sum \frac{r_i}{x-\gamma_i}, \ m_i, r_i$ integers.    Let $\tilde{f} = e^{\int h^*}$ and $\tilde{g} = e^{\int g}$.  We have that $f = \tilde{f} \frac{\sigma(\tilde{g})}{\tilde{g}}$.  Since the polar dispersion of $h^*$ is $0$, we have that $\tilde{f}$ is standard. \\[0.1in]
Assume $\sigma(q) = qx$.  Let $\dd = x\frac{d}{dx}$ and
\[ h = \frac{\dd(f)}{f} = \sum_i\frac{n_ix}{x-\alpha_i}\]
for some integers $n_i$.  Note that
\[ \frac{nx}{x-q^{-N}\alpha} = (\sigma(\frac{nx}{x-q^{-(N-1)}\alpha}) -  \frac{nx}{x-q^{-(N-1)}\alpha}) + \frac{nx}{x-q^{-(N-1)}\alpha}.\]
We can therefore proceed as in the second part of 1.~above and find that $h = h^* + \sigma(g) - g$, where $\pdisp(h^*) = 0$ and $h^* = \sum \frac{m_ix}{x-\beta_i}, \ g = \sum \frac{r_ix}{x-\gamma_i}, \ m_i, r_i$ integers.    Let $\tilde{f} = e^{\int h^*/x}$ and $\tilde{g} = e^{\int g/x}$.  We have that $f = \tilde{f} \frac{\sigma(\tilde{g})}{\tilde{g}}$.  Since the polar dispersion of $h^*$ is $0$, we have that $\tilde{f}$ is standard. \end{proof}

\begin{lem}\label{posdisp} Let $k$ be as above and $f \in k , a \in C, a\neq 0$. \begin{enumerate}
\item If $\sigma(x) = x+1$,  and $f$ has a pole, then $\pdisp(\sigma(f) - af) >0$.
\item If $\sigma(x) = qx$ and $f$ has a nonzero pole, then $\pdisp(\sigma(f) - af) >0$.
\end{enumerate}\end{lem}
\begin{proof} We will prove case 2.; the proof of case 1.~is similar (a proof of this latter case, when $a$ = 1,  also appears in \cite{matusevich}).  Let $\pdisp(f)= N \geq 0$ and let $\alpha$ be a  pole of $f$ such that $q^N\alpha$ is also a pole of $f$.  One sees that $q^{-1}\alpha$ is a pole of $\sigma(f)$ but not of $f$ and $q^{N}\alpha$ is a pole of $f$ but not of $\sigma(f)$.  Therefore $\pdisp(\sigma(f) - f) \geq N+1 > 0$  \end{proof}

\begin{lem} \label{trivlem}Let $k=C(x)$, $x$ transcendental over $C$, $0\neq a \in C$. \begin{enumerate}
\item Assume $\sigma(x) = x+1$  and $\dd = \frac{d}{dx}$.  Let $b\in k$ and assume that there is a nonzero $L \in C[\dd]$ such that $L(b) = h(x+1) - ah(x)$ for some $h \in k$. Then $b = \sigma(e) - ae$ for some $e \in k$.
\item Assume $\sigma(x) = qx$  and $\dd = x\frac{d}{dx}$.  Let $b \in k$ and assume that there is a nonzero $L \in C[\dd]$ such that $L(b) = h(qx) - ah(x) $ for some $h \in k$. Then
\begin{enumerate}
\item if $a = q^r$ for some $r \in \ZX$, then $b = \sigma(e) - ae +cx^r$ for some $e \in k$ and $c \in C$, or, if not,
\item  $b = \sigma(e) - ae $ for some $e \in k$.
\end{enumerate}
\end{enumerate}
\end{lem}

\begin{proof} 1.~Let $b = b^*+\sigma(g) - ag$ as in Lemma~\ref{normalform}.  We then have that $L(b^*) = \sigma(h- L(g))
-a(h-L(g))$. We shall first show that $b^*$ must be a polynomial.  Since $\pdisp(b^*) = 0$ and all the poles of $L(b^*)$ appear as poles of $b^*$, we must have that $\pdisp(L(b^*)) = 0$.  If $b^*$ has a pole, then so must $L(b^*)$.  This furthermore implies that $h-L(g)$ also has a pole.  This yields a contradiction since, in this case, Lemma~\ref{posdisp} implies $\pdisp(\sigma(h- L(g)) - a(h - L(g)) >0$. \\[0.1in]
Let $C[x]^{\leq n}$ denote the vector space of polynomials of degree at most $n$. If $a \neq 1$, then the map $p \mapsto \sigma(p) - ap$ from $C[x]^{\leq n}$ to $C[x]^{\leq n}$ has trivial kernel and so must be surjective.  If $a = 1$, this map has kernel $C$ and maps  surjectively onto $C[x]^{\leq n-1}$.  In either case, this implies that for any polynomial $p$, there is a polynomial $q$ such that $\sigma(q) - aq = p$. Therefore, we may write $b^* = \sigma(\bar{h}) - \bar{h}$ for some polynomial $\bar{h}$ and so $b = \sigma(g+\bar{h}) - (g+\bar{h})$.\\[0.1in]
\noindent 2.~Let $b = b^*+\sigma(g) - g$ as in Lemma~\ref{normalform}.  We then have that $L(b^*) = \sigma(h- L(g))
-a(h-L(g))$. We shall first show that $b^* \in \CX[x, \frac{1}{x}]$.  Since $\pdisp(b^*) = 0$ and all the poles of $L(b^*)$ appear as poles of $b^*$, we must have that $\pdisp(L(b^*)) = 0$. If $b^*$ has a nonzero pole, then so must $L(b^*)$. This furthermore implies that  $h-L(g)$ also has a nonzero pole. This yields a contradiction since, in this case, Lemma~\ref{posdisp} implies $\pdisp(\sigma(h- L(g)) - (h - L(g)) >0$. \\[0.1in]
For any integer $r$, we have that $\sigma(x^r) - ax^r = (q^r - a) x^r$.  If $a = q^r$, then  any $\bar{g} \in \CX[x, \frac{1}{x}]$ without a term of the form $cx^r$ may be written as $\bar{g} = \sigma(\bar{h}) - a\bar{h}$ for some $\bar{h} \in \CX[x,\frac{1}{x}]$.  In particular, there exists a $c \in \CX$ such that $b^* -cx^r = \sigma(\bar{h}) - a\bar{h}$ and so $b = \sigma(g+\bar{h}) - a(g+\bar{h}) + cx^r$.  Conversely, if $b = \sigma(e) - e + cx^r$, then $\dd(b) - r b = \sigma(\dd(e) - re)  - (\dd(e) - r e)$.  If $a \neq q^r$, then any $\bar{g} \in \CX[x, \frac{1}{x}]$  may be written as $\bar{g} = \sigma(\bar{h}) - a\bar{h}$ for some $\bar{h} \in \CX[x,\frac{1}{x}]$.
In particular,  $b^*  = \sigma(\bar{h}) - a\bar{h}$ and so $b = \sigma(g+\bar{h}) - a(g+\bar{h}) $. \end{proof}

One more technical result is needed in the text.
\begin{lem}\label{descent} Let $k \subset k(u)$ be $\sigma$-fields with $u$ transcendental over $k$ and $\sigma(u) = au$ for some $a \in k$.  Let $h, b \in k$ and assume that $\sigma(y) - by = h$ has a solution in $k(u)$.  Then $\sigma(y) - by = h$ has a solution in $k$.
\end{lem}
\begin{proof} Let $ w \in k(u)$ satisfy $\sigma(w) - pw = h$. We may write $w = \frac{p}{q} + r$ where $p,q,r \in k[u]$ and $\deg_uq <\deg_p$. Note that $\frac{\sigma(p)}{\sigma(q)} - b \frac{p}{q} = \frac{\tilde{p}}{\tilde{q}}$ where $\deg_u\tilde{p}<\deg_u\tilde{q}$.  Therefore $h$ must equal the $u^0$ term in $\sigma(r) - br$ which is $\sigma(r_0) - b r_0$ where $r_0$ is the $u^0$ term in $r$.\end{proof}

\subsection{$\SDP$-Galois Theory}\label{appendix} In this section we develop a Galois theory of systems of differential-difference equations that measures the differential properties  (with respect to an auxiliary set of derivations) of solutions of this system.  We will show how the results of Section~\ref{galoistheory} (as well as the usual difference and differential Galois theories and the parameterized differential Galois theory of \cite{CaSi}) are special cases of this general theory.  This general theory also yields a Galois theory of completely integrable differential/difference equations, a result that seems also to be new. Many of the statements and proofs follow along the lines of the usual Picard-Vessiot theories of differential or difference equations \cite{PuSi, PuSi2003} but in this theory there are some new complications (especially in the existence of Picard-Vessiot extensions). We therefore present the theory {\em ab initio} and give complete proofs for the important results. For another approach to combining difference and differential Galois theories, see \cite{andre_galois}. We begin with some definitions that generalize the corresponding definitions in these cases.

\begin{defins}\label{def6.1} (a) A {\em $\SDP$-ring}  is a ring $R$ with a set of automorphisms $\Sigma = \{\sigma_i\}$ and sets of derivations $\Delta = \{\delta_i\}$ and $\Pi = \{\der_i\}$ such that for any $\tau, \mu \in \Sigma\cup\Delta\cup\Pi, \ \mu(\tau(r)) = \tau (\mu(r)) \ \forall \ r \in R$. The notions of $\SDP$-field, $\SDP$-ideal, $\SDP$-homomorphism, etc. are defined similarly. All fields considered will be of characteristic zero.

(b)  The {\em $\SDP$-constants $R^\SDP$} of a $\SDP$-ring $R$  is the set  $R^\SDP = \{ c \in R \ |  \delta(c) = \der(c) = 0 \ \forall \delta\in \Delta, \der \in \Pi \mbox{ and } \sigma(c) = c \ \forall \sigma \in \Sigma\}$. Similarly, the {\em $\SD$-constants} are  $R^\SD = \{ c \in R \ |  \delta(c)  = 0 \ \forall \delta\in \Delta,  \mbox{ and } \sigma(c) = c \ \forall \sigma \in \Sigma\}$.

(c) A {\em simple $\SDP$-ring}  is a $\SDP$-ring whose only $\SDP$-ideals are $(0)$ and $R$.

 (d) Given a $\SDP$-field $k$, a {\em $\SD$-linear system} is  a system of equations of the form
\begin{eqnarray*}
\sigma_i(Y) &=& A_iY \ A_i \in \GL_n(k),  \ \sigma_i \in \Sigma\\
\delta_i(Y) &=& B_iY \ B_i  \in \gl_n(k), \ \delta_i \in \Delta
\end{eqnarray*}
where the $A_i, B_j$ satisfy the integrability conditions
\begin{eqnarray*}
\sigma_i(A_j)A_i &=& \sigma_j(A_i)A_j\\
\sigma_i(B_j) A_i & = & \delta_j(A_i) + A_iB_j\\
\delta_i(A_j) + A_jA_i &= &\delta_j(A_i) + A_i A_j
\end{eqnarray*}
for all $\sigma_i, \sigma_j \in \Sigma$ and all $\delta_i,\delta_j \in \Delta$.
\end{defins}

We note that the integrability conditions are precisely the
conditions that must be met if $\mu(\tau(Z)) = \tau(\mu(Z))$ for any
solution $Z$ of the system (in some $\SDP$-extension ring) and any
$\mu,\tau \in \Sigma \cup \Delta$. The above definitions are broad
enough to include many of the linear systems considered by Galois
theories.  For example when
\begin{itemize}
\item[-] $\Sigma$ and $\Pi$ are empty, one is just  considering differential fields and an integrable system of linear differential equations as in \cite{PuSi2003}, Appendix D.
 \item[-] $\Delta$ and $\Pi$ are empty, one is considering difference fields and systems of ``integrable'' difference equations, a generalization of the ordinary linear difference equations considered in \cite{PuSi}.
 \item[-] $\Pi$ is empty, one is considering the differential/difference fields  and  linear functional equations of  \cite{wu05} and \cite{blw}.
\item[-] $\Sigma$ is empty, one is considering the parameterized  differential equations of \cite{CaSi}.
\item[-]  $\Sigma$ and $\Pi $ are singletons, one is considering the $\sd$ situation considered in the first part of this paper.
\item[-] $\Pi$ is empty, one is considering integrable differential/difference linear systems.
\end{itemize}

The above observation shows that the usual examples of linear differential equations, linear difference, parameterized differential equations give examples of $\SD$-linear systems over $\SDP$-fields.   We have the following simple facts:

\begin{lem}\label{radical} Let $R$ be a $\SDP$-ring.\begin{enumerate}
\item The radical of a $\SDP$-ideal  is a $\SDP$-ideal.
\item A maximal $\SDP$-ideal $I$ is radical and furthermore for any $r\in R$, $r \in I$ if and only if $\sigma(r) \in I \ \forall \sigma \in \Sigma$.  Therefore, $R/I$ is a reduced $SDP$-ring. \end{enumerate}
\end{lem}
\begin{proof} 1.~follows from the corresponding statements for difference ideals (obvious)  and differential ideals (Lemma 1.8, \cite{kaplansky}).\\
2.~follows from (a) and the fact that $\{r \in R \ | \sigma(r) \in I \ \forall \sigma \in \Sigma \}$ is a $\SDP$-ideal containing $I$.  \end{proof}

\begin{lem}\label{idempotents} Let $k$ be a $\SDP$-field and $R$ a simple $\SDP$-ring, finitely generated over $k$ as a $\Delta\Pi$-ring.  Then there exist idempotents $e_0, \ldots ,e_{t-1}$ such that \begin{enumerate}
\item $R = R_0\oplus \ldots \oplus R_{t-1}, \ R_i = e_iR$,
\item each $\sigma \in \Sigma$ permutes the set $\{R_0, \ldots ,R_{t-1}\}$ and $<\Sigma>$ acts transitively on this set, where $<\Sigma>$ is the semigroup generated by $\Sigma$.  Furthermore  $\sigma^t$ leaves each $R_i$ invariant, and
\item each $R_i$ is a domain and is a simple $\tilde{\Sigma}\Delta\Pi$-ring where $\tilde{\Sigma} = \{\sigma^t \ | \ \sigma \in \Sigma\}$
\end{enumerate}
\end{lem}
\begin{proof} This lemma and its proof are straightforward generalizations of Corollary 1.16 of \cite{PuSi}.  Since $R$ is reduced and finitely generated as a differential $k$-algebra, $(0)$ is a radical ideal and the Ritt-Raudenbush Theorem (\cite{DAAG}, p.126) allows us to write $(0) = \cap_{i=0}^{t-1}I_i$, where each $I_i$ is a prime $\Delta\Pi$-ideal. Assuming this representation is irredundant, it is unique as well.  For any $\sigma \in \Sigma$, $(0) =  \cap_{i=0}^{t-1}\sigma(I_i)$, so $\sigma$ permutes the $I_i$.  Furthermore, $<\Sigma>$ must act transitively on $\{I_0, \ldots, I_{t-1}\}$ or else the intersection of the elements of a proper orbit would be a nontrivial $\SDP$-ideal.  Since the elements of $\Sigma$ commute, $<\Sigma>$ acts transitively on the orbits of any $\sigma \in \Sigma$.  Therefore for any $\sigma \in \Sigma$, the orbits of $\sigma$ have the same size, which must therefore divide $t$. We can conclude therefore that $\sigma^t$ leaves each $I_i$ invariant.

Fix some $\sigma \in \Sigma$ and, for each $i$, let $J_i = \{r \in R \ | \ \sigma^t(r) \in I_i\}$.
For each $i$, $J_i$ is a prime $\Delta\Pi$-ideal containing $I_i$.  Any $\tau \in \Sigma$ permutes the $J_i$ so $\cap J_i$ is a proper $\SDP$-ideal of $R$ and therefore equals $(0)$.  Therefore for each $i, J_i = I_i$ so $r \in I_i$ if and only $\sigma^t(r) \in I_i$.  This implies that the ring $S_i = R/I_i$ is a $\tilde{\Sigma}\Delta\Pi$-integral domain where $\tilde{\Sigma} = \{\sigma^t \ | \ \sigma \in \Sigma\}$.

We now claim that for each $i$, $S_i$ has no nonzero proper $\tilde{\Sigma}\Delta\Pi$-ideals. It is enough to show that $R$ has no proper $\tilde{\Sigma}\Delta\Pi$-ideals properly containing $I_i$. Assume the contrary and let $J_i$ be such an ideal.  Let $N_0 = J_i, N_1, \ldots , N_r$ be the orbit of $J_i$ under $<\Sigma>$.  Note that since $J_i$ is left invariant under $\sigma^t$ for all $\sigma \in \Sigma$, we know that this orbit is finite.  Since $\cap N_i$ is a proper $\SDP$-ideal it must be $(0) \subset I_i$.  Since $I_i$ is prime, we must have that some $N_\ell \subset I_i$.  We may write $N_\ell = \theta(J_i)$ where $\theta$ is a powerproduct of the $\sigma \in \Sigma$.  Therefore $\theta(I_i) \subset \theta(J_i)  \subset I_i$. Since $\theta(I_i)$ must be one of the $I_i$, and these are not contained in each other, we have $\theta(I_i) = I_i$.  Therefore, $J_i \subset \theta^{-1}(I_i) = I_i$, a contradiction.

We now claim that the $I_j$ are pairwise comaximal, that is, $I_i+I_j = R$ if $i \neq j$.  To see this note that $I_i+I_j$ is a $\tilde{\Sigma}\Delta\Pi$-ideal properly containing $I_i$ and $I_j$.  By the above, it must be all of $R$.  Let $p_i:R \rightarrow  S_i$ be the canonical projection. The Chinese Remainder Theorem (\cite{zariski_samuel}, Ch. 3, Sec. 13) implies that the map $p:R \rightarrow \sum_{i = 0}^{t-1} S_i$ given by $p(r) = (p_0(r), \ldots ,p_{t-1}(r))$  is a ring homomorphism.  We may write $R = \sum_{i=0}^{t-1} R_i$ where $R_i = p^{-1}(S_i)$.  \end{proof}

In Section~\ref{galoissec}, we will state and prove a Galois correspondence.  For this we will need the following corollary.

\begin{cor}\label{zerodivisor} Let $k, R, R_i$ and $e_i$ be as in Lemma~\ref{idempotents} and $K$ the total ring of quotients of $R$.  Then
\begin{itemize}
\item $K = K_0\oplus \ldots \oplus K_{t-1}$, where $K_i$ is the quotient field of $R_i$.
\item If $S$ is an $\Sigma$-subring of $K$, then any $s\in S$ that is a zero divisor in $K$ is a zero divisor in $S$. In particular, we can embed the total ring of quotients of $S$ into $K$.
\item If $S$ has the property that any non-zero divisor in $S$ is invertible in $S$, then there exists a partition $X_0 \cup \ldots \cup X_{d-1} = \{e_0, \ldots e_t\}$ such that for $E_i = \sum_{e_j \in X_j} e_j$
\begin{enumerate}
\item $S = \oplus_{i=1}^{d-1} SE_i$ where each $SE_i$ is a field, and
\item the $SE_i$ are left invariant  by each $\sigma^d, \ \sigma \in \Sigma$,
\end{enumerate}
\end{itemize}\end{cor}

\begin{proof} ({\em cf.,} \cite{PuSi}, proof of Theorem 1.29) Item 1.~is clear from the description of $R$ in the previous lemma.  To prove items 2.~and 3., we make the following definition. A subset $A \subset \{0, \ldots , t-1\}$ is called a support for $S$ if there is an element $f \in (f_0, \ldots , f_{t-1}) \in S$ such that $A = \{i \ | f_i \neq 0\}$. We shall write $A = supp(f)$ in this case.  If $A$ and $B$ are supports then so are $A\cap B$ and $A\cup B$. The semigroup $<\Sigma>$ permutes the $R_i$ and so acts on $\{0, \ldots , t-1\}$.  One sees that if $A$ is a support and $\sigma \in <\Sigma>$, then $\sigma(A)$ is a support as well.

Let $A$ be a minimal support containing $0$ and assume that $A$ is the support of $f \in S$.  For any $\phi_1, \phi_2 \in <\Sigma>$, the minimality of $A$ implies that $\phi_1(A) = \phi_2(A)$ or $\phi_1(A) \cap \phi_2(A) = \emptyset$.  Since $<\Sigma>$ acts transitively on $\{0, \ldots , t-1\}$, we have that the disjoint sets $X_0, \ldots , X_{d-1}$ in $\{\phi(A)\ | \ \phi \in <\Sigma>\}$ form a partition of $\{0, \ldots , t-1\}$ by sets each of the same cardinality $t/d$ (which is the cardinality of $A$).  Note that each $X_i$ is also a support in $S$.

Let $g \in S$ be a zero divisor in $K$.  This implies that $supp(g) \neq \{0, \ldots , t-1\}$.  Furthermore, we must have $supp(g) \cap X_i$ is a proper subset of $X_i$ for some $i$.  By minimality of $A$, we must have that $supp(g) \cap X_i = \emptyset$.  If $X_i = supp(h)$ for some $h \in S$, then $gh = 0$, so $g$ is a zero divisor in $S$.  This proves statement 2.

Now assume that any non-zero divisor in $S$ is invertible in $S$.  Let $f \in S$ have minimal support $A$.  As above, there are $\phi_0 = id, \phi_1, \ldots ,\phi_{d-1} \in <\Sigma>$ such that the $\{\phi_i(f)\}$ have disjoint supports whose union is all of $\{0, \ldots , t-1\}$.  Therefore $r =  f+\phi_1(f)+ \ldots +\phi_{t-1}(f)$ is a nonzero divisor and so is invertible.  Letting $E_0 = r^{-1}f$ we have $E_0 = \sum_{j \in A} e_j$.  Similarly we define $E_i = r^{-1}\phi_i(f)$.  These $E_i$ satisfy the properties of statement 3. \end{proof}

\subsubsection{$\Sigma\Delta\Pi$-Picard-Vessiot Extensions} As in the usual Picard-Vessiot theory of differential or difference equations, we shall consider ``splitting rings'' for our equations.
 Let $k$ be a $\SDP$-field and let
\begin{eqnarray}
\sigma_i(Y) &=& A_iY \ A_i \in \GL_n(k),  \ \sigma_i \in \Sigma \label{system}\\
\delta_i(Y) &=& B_iY \ B_i  \in \gl_n(k), \ \delta_i \in \Delta \notag
\end{eqnarray}
be a $\SD$-linear system over $k$.
\begin{defin}\label{defPV} A $\SDP$-ring $R$ is a {\em $\SDP$-Picard-Vessiot ring ($\SDP$-PV ring)} for equations (\ref{system}) if
\begin{enumerate}
\item $R$ is a simple $\SDP$-ring,
\item there exists a $Z\in \GL_n(R)$ such that $\sigma_i(Z) = A_iZ \ \forall\sigma_i \in \Sigma $ and $\delta_i(Z) = B_iZ \ \forall \delta_i \in \Delta$, and
\item $R = k\{Z, \frac{1}{\det A}\}_\Pi$, that is $R$ is generated a  $\Pi$-ring by the entries of $Z$ and the inverse of the determinant of $Z$.
\end{enumerate}
\end{defin}

Given a $\SDP$-field $k$ and a linear system (\ref{system}), one can show that a $\SDP$-PV ring for this system always exists. To do this, let $Y = (y_{i,j})$ be an $n\times n$ matrix of $\Pi$-differential indeterminates and let $S = k\{Y,\frac{1}{\det Y}\}_\Pi$.  We extend the action of $\Sigma$ and $\Delta$ to this ring by prolonging the equations of (\ref{system}) using the derivations of $\Pi$.  For example, for any $\sigma_i \in \Sigma$ and $\der \in \Pi$ we let $\sigma_i(Y) = A_i Y, \sigma_i(\der Y ) = (\der A_i)Y + A_i(\der Y), \ldots $ and for any $\delta_i \in \Delta$, $\delta_i(Y) = B_i Y, \delta_i(\der Y) = (\der B_i)Y + B_i (\der Y), \ldots $.  In this way $S$ becomes a $\SDP$-ring.  Letting $I$ be a maximal $\SDP$-ideal, $R=S/I$ is a $\SDP$-PV-ring.

In the usual Galois theory of linear difference or differential equations, one needs to assume that the constants of the ground field are algebraically closed before one can assert that  Picard-Vessiot extensions  are unique or that there are enough automorphisms to yield a Galois correspondence.  In the Parameterized Picard-Vessiot theory of $\cite{CaSi}$, one needs the addition hypothesis that the appropriate field of constants is differentially closed (with respect to the parametric derivatives).  A similar condition here guarantees  uniqueness and existence of enough automorphisms.  Before we formally state and prove this we need some preliminary lemmas.

\begin{lem} \label{disjoint} If $k \subset K$ are $\SD$-fields then $k$ and $K^\SD$ are linearly disjoint over $k^\SD$.
\end{lem}
\begin{proof} It is enough to show that if $c_1, \ldots , c_n \in K^\SD$ are linearly dependent over $k$, then they are linearly dependent over $k^\SD$.  Let $\sum_{i=1}^n a_ic_i = 0, \ a_i \in k$. Among all such relations, select one with the minimal number of nonzero $a_i$. We can assume that $a_1 = 1$ and will show that this implies that the other $a_i$ are in $k^\SD$. For any $\sigma \in \Sigma, \ \delta \in \Delta$, we have $\sum_{i=2}^n (\sigma(a_i) - a_i)c_i = 0$ and $\sum_{i=2}^n \delta(a_i)c_i =0$, so by minimality $a_i \in k^\SD$.\end{proof}

\begin{lem}\label{ideals1} Let $k$ be a $\SDP$-field and $R$  be a simple $\SDP$-ring containing $k$ that is finitely generated as a $\Pi$-ring over $k$.  Let $K$ be the total ring of quotients of $R$ and let $K\{y_1, \ldots ,y_n\}_\Pi$ be the ring of $\Pi$-differential polynomials over $K$. Extend the action of $\Sigma$ and $\Delta$ to $K\{y_1, \ldots ,y_n\}_\Pi$ by setting $\sigma(y_i) = y_i  \ \forall \sigma \in \Sigma$ and $\delta(y_i) = 0 \ \forall \delta \in \Delta$.  Then the map  $I \mapsto(I)$ is a bijective correspondence from the set of $\Pi$-ideals of $K^\SD\{y_1, \ldots ,y_n\}_\Pi$ and the  $\SDP$-ideals of $K\{y_1, \ldots ,y_n\}_\Pi$. In particular, letting $R=k$, we have that the map  $I \mapsto(I)$ is a bijective correspondence from the set of $\Pi$-ideals of $k^\SD\{y_1, \ldots ,y_n\}_\Pi$ and the  $\SDP$-ideals of $k\{y_1, \ldots ,y_n\}_\Pi$.
\end{lem}
\begin{proof} ({\em cf.}, \cite{PuSi2003}, Lemma 1.23) The hypothesis that $R$ is a simple $\SDP$-ring containing $k$ that is finitely generated as a $\Pi$-ring over $k$ is only used to conclude that $K = \oplus_{i=0}^{t-1}e_i K$ where the $e_i$ are a maximal set of orthogonal idempotents,  the $ K_i = e_iK$ are $\Delta\Pi$-fields and $<\Sigma>$ acts transitively on the set $\{K_0, \ldots , K_{t-1}\}$.  Therefore the lemma could be stated in greater generality but this will not be needed.

It is enough to show that any $\SDP$-ideal $J$ is generated by $J \cap k^\SD\{y_1, \ldots ,y_n\}_\Pi = I$.  Any   $f \in J$ can be written  as $f = \sum f_\alpha y_\alpha$ where the $f_\alpha \in k$ and the $y_\alpha$ are differential monomials (\cite{DAAG}, p. 70), that is, power products of derivatives of the $y_i$. Define $\ell(f)$ to be the number of $f_\alpha \neq 0$.  We will proceed, by induction on $\ell(f)$, to show that $f \in (I)$. If $\ell(f) = 0,1$, this is obvious.  Assume $\ell(f) > 1$ and that the claim is true for elements of smaller length.

We will first show that $f \in (I)$ under an additional
hypothesis on the form of $f$.  We will assume that there is a  $j, 0\leq j \leq t-1$ such that  $f = \sum \theta_i(e_jh)$ where   $h \in K\{y_1, \ldots ,y_n\}_\Pi, \ e_jh \neq 0$ and each $\theta_i \in <\Sigma>$ has the property that $\theta_i(e_j) = e_i$ and $\theta_j = id$. Note that $e_jh \neq 0$ implies that the nonzero coefficients of this polynomial are invertible in $K_j$. Therefore any nonzero coefficient of $f$ is invertible, so we may assume  that some $f_\alpha = 1$. If all the $f_\alpha \in K^\SD$ we are done so assume some $f_{\alpha_1} \in K\backslash K^\SD$. Assume $\sigma(f_{\alpha_1}) \neq f_{\alpha_1}$ for some $\sigma \in \Sigma$ (the argument if $\delta(f_{\alpha_1}) \neq f_{\alpha_1}$ for some $\delta \in \Delta$ is similar). We have that $\ell(\sigma(f) - f) < \ell(f)$ and that $\ell(\sigma(f_{\alpha_1}^{-1}f) - f_{\alpha_1}^{-1}f) < \ell(f)$. Since $\sigma(f_{\alpha_1}^{-1}f) - f_{\alpha_1}^{-1}f = \sigma(f_{\alpha_1}^{-1})(\sigma(f)-f) + (\sigma(f_{\alpha_1}^{-1})-f_{\alpha_1}^{-1})f$, we have $(\sigma(f_{\alpha_1}^{-1})-f_{\alpha_1}^{-1})f \in (I)$. Note that $(\sigma(f_{\alpha_1}^{-1})-f_{\alpha_1}^{-1})$ may be a zero divisor so we cannot immediately conclude that $f\in (I)$.  Nonetheless, we have that for some $i, \ e_i f \in (I)$. Note that $(I)$ is not only invariant under each $\sigma \in \Sigma$ but also under each $\sigma^{-1}$. Therefore, for each $l, \ \theta_l \theta_i^{-1}(e_if) = \theta_l \theta_i^{-1}(\theta_i(e_jh) = \theta_l(e_jh) = e_lf$. This implies that $e_lf \in (I)$ for all $l$ so $f \in (I)$.

We now consider general $h \in J$.  To show $h \in (I)$ it is enough to show that each $e_jh \in (I)$.  Fix a $j$. For each $j$ there exists a $\theta_i \in <\Sigma>$ such that $\theta_i(e_j) = e_i$.  Let $f = \sum_{i=0}^{t-1}\theta_i(e_jh)$.   We have shown that $f \in (I)$ so  $e_j f = e_jg\in (I)$.  Proceeding in this way for all $j$, we can then conclude that $f \in (I)$. \end{proof}

Finally, we will need the following result  (which is a special case of a result of Kolchin (\cite{DAAG}, Theorem 3, p. 140) concerning the extension of homomorphisms; for those familiar with model theory, this result can be replaced by an appeal to elimination of quantifiers for differential fields.  Let $k$ be a $\Pi$-differential field for some set of commuting  derivations $\Pi$. A {\em semiuniversal extension $U$ of $k$} is an $\Pi$-differential field $U$ containing $k$ such that any finitely generated $\Pi$-extension of $k$ can be embedded into $U$. Kolchin has shown that such extensions exist (\cite{DAAG}, p.92).

\begin{prop} \label{kolchin} ( {\em cf.}~\cite{DAAG}, Theorem 3, p. 140) Let $k$ be a $\Pi$-differential field, $R$ a $\Pi$-integral domain, $\Pi$-finitely generated over $k$, and $R_0$ a $\Pi$-ring with $k \subset R_0 \subset R$.  There exists a nonzero element $u_0$ of $R_0$ such that every differential $k$-homomorphism $\phi$ of $R_0$ into $U$, a semiuniversal extension of $k$, with $\phi(u_0) \neq 0$ can be extended to a differential k-homomorphism $\tilde{\phi}$ of $R$ into $U$. \end{prop}
The following will be used to show  that $\SDP$-PV extensions of fields with $k^\SD$ differentially closed  are unique and that we will have enough automorphisms.

\begin{prop}\label{nonewconstants} Let $k$ be a $\SDP$-field with $k^\SD$ differentially closed as a $\Pi$-field and let $R$ be a simple $\SDP$-ring containing $k$ that is finitely generated as a $\Pi$-ring over $k$.  Then $R^\SD = k^\SD$.
\end{prop}
\begin{proof}We first prove this proposition when $R$ is an integral domain. We will argue by contradiction. Let $c \in R^\SD\backslash k^\SD$.  We will show that there is a $\Pi$-$k$-algebra homomorphism $\tilde{\phi}$ of $R$ into a semiuniversal $\Pi$-extension $U$ of $k$ such that  $\tilde{\phi}(c)  = \tilde{c}\in k^\SD$.  Assuming this for a moment, we will show how we get a contradiction.

Since $c, \tilde{c} \in  R^\SD$, the $\Pi$-ideal $(c-\tilde{c})_\Pi$ is a nonzero $\SDP$-ideal as well.  Since $R$ is simple, we must have $1 \in (c-\tilde{c})_\Pi$.  Since $\tilde{\phi}(c-\tilde{c})= 0$ and $\tilde{\phi}$ is a $\Pi$-homomorphism, we have a contradiction.

To show the existence of $\tilde{\phi}$, let $R_0 = k\{c\}_\Pi \subset R$. Note that $R_0$ is again a $\SDP$-ring.  Applying Proposition~\ref{kolchin} to $R_0 \subset R$ considered as $\Pi$-rings, there is a nonzero $u_0 \in R_0 $  satisfying the conclusions of this proposition. We may write $u_0 = U_0(c)$ for some $U_0 \in k\{y\}_\Pi$, the ring of $\Pi$-differential polynomials.  Let $J \subset k\{y\}_\Pi$ be the defining $\Pi$-ideal of $c$ over $k$.  Note that since $c \in R^\SD$, $J$ is a $\SDP$-ideal of $k\{y\}_\Pi$ when we extend $\Sigma$ and $\Delta$ to $k\{y\}_\Pi$ as in Lemma~\ref{ideals1}. This lemma implies that $J$ is generated by $I = J \cap k^\SD\{y\}_\Pi$.  Let $U_0 = \sum a_i V_I$ where $\{a_i\}$ is a $k^\SD$-basis of $k$ over $k^\SD$ and the $V_i \in k^\SD\{y\}_\Pi$.  Since $U_0(c)\neq 0$ we have that $V_i(c) \neq 0$ for some $i$.  Since $k^\SD$ is a $\Pi$-differentially closed field and there is a zero of $I$ such that $V_i \neq 0$ in some extension field, we have that there is a $\tilde{c} \in k^\SD$ that is a zero of $I$ and such that $V_i(\tilde{c}) \neq 0$. Lemma~\ref{disjoint} implies that $U_0(\tilde{c}) \neq 0$.  Therefore the map $c\mapsto \tilde{c}$ defines a $\Pi$-homomorphism $\phi$ of $R_0$ to $k$ with $\phi(c) \in k^\SD$ and $\phi(u_0) \neq 0$.  Proposition~\ref{kolchin} implies that we can extend $\phi$ to a map $\tilde{\phi}$ yielding the desired conclusion.

We now remove the hypothesis that $R$ is an integral domain.  We may write $R = R_0\oplus \ldots \oplus R_{t-1}$ as in Lemma~\ref{idempotents}, where $R_i = e_iR$ and the $e_i$ form a set of mutually orthogonal idempotents with $1 = e_0+\ldots +e_{t-1}$.  We note that this implies that any $\sigma \in \Sigma$ permutes the $e_i$ and that, since $\sigma^t$ leaves each $R_i$ invariant, $\sigma^t(e_i) = e_i$.  Let $r \in R^\SD$ and write $r = r_0+\ldots + r_{t-1}$ where $r_i = e_i r$. We then have, for any $\sigma \in \Sigma$, that $\sigma^t(r) = r$ so $\sigma(r_i) = r_i$.  Lemma~\ref{idempotents} implies that each $R_i$ is a simple $\tilde{\Sigma}\Delta\Pi$-{\em domain} ($\tilde{\Sigma} = \{ \sigma^t \ | \ \sigma\in \Sigma$), so the above argument implies that  $r_i \in e_ik^{\tilde{\Sigma}\Delta}$.  Since $k^\SD$ is differentially closed it must be algebraically closed and we can conclude that $k^\SD = k^{\tilde{\Sigma}\Delta}$.  Therefore $r = e_0s_0 + \ldots +e_{t-1}s_{t-1}$ for some $s_i \in k^\SD$.  Since $\Sigma$ acts transitively on the $e_i$, we have  that $s_0 = s_1 = \ldots =s_{t-1}$ so $s \in k^\SD$.\end{proof}
\begin{cor} \label{cornonewconstants} Let $k$ be a $\SDP$-field with $k^\SD$ differentially closed as a $\Pi$-field and let $R$ be a simple $\SDP$-ring containing $k$ that is finitely generated as a $\Pi$-ring over $k$.  If $K$ is the total quotient ring of $R$, then $K^\SD = k^\SD$.
\end{cor}
\begin{proof} Let $c = \frac{a}{b} \in K^\SD$.  Using Proposition \ref{nonewconstants}, it suffices to show that $R\{c\}_\Pi$ is a simple $\SDP$-ring, since this would imply that $c \in k^\SD$.  Let $J$ be a nonzero $\SDP$-ideal of $R\{c\}_\Pi$.  We claim that $J \cap R$ contains a nonzero element.   Assuming that this is the case, then, since $R$ is simple, $J\cap R = R$ so $J =R\{c\}_\Pi$.    To prove the claim let $0 \neq u \in J$.  We may write $u = \sum a_i \theta_i (c)$ where each $\theta_i$ is a power product of elements in $\Pi$.  Using elementary properties of derivations, one sees that for each $i$, there is a positive integer $n_i$ such that $b^{n_i}\theta_i(c) \in R$.  Therefore there is a positive integer $n$ such that $b^n u \in R$.  Since $b$ is not a zero divisor,  $b^n u $ is a nonzero element of $J\cap R$. \end{proof}

From Proposition~\ref{nonewconstants}, we can conclude the uniqueness of $\SDP$-PV extensions over a $SDP$-field with $\Pi$-differentially closed $k^\SD$.
\begin{prop} \label{unique} Let $k$ be a $\SDP$-field with $k^\SD$ a $\Pi$-differentially closed field.  Let $R_1$ and $R_2$ be $\SDP$-PV extensions of $k$ for the linear system (\ref{system}).  Then there exists a $\SDP$-k-isomorphism between $R_1$ and $R_2$

\end{prop}
\begin{proof}The proof is the same as the proof of  Proposition 1.9 on page 7 of \cite{PuSi}.
\end{proof}

We end this section by  characterizing the total ring of quotients of a $\SDP$-PV extension $R$ of a $\SDP$-field $k$.  We say that a $\SDP$-ring $K$ containing $k$ is a {\em total $\SDP$-PV ring} if it is the total ring of quotients of a $\SDP$-PV ring over $k$

\begin{prop}\label{totalquotprop} Let $k$ be a $\SDP$-field with $k^\SD$ a $\Pi$-differentially closed differential field.  Let $K\supset k$ be a $\SDP$ring satisfying
\begin{enumerate}
\item $K$ has no nilpotent elements and every non-zero divisor of $K$ is invertible.
\item $K^\SD = k^\SD$.
\item $K = k<Z>_\Pi$ where $Z \in \GL_n(K)$ satisfies equations (\ref{system})

\end{enumerate}
Then $R$  is a $\SDP$-PV ring and $K$ is a total $\SDP$-PV ring for equations (\ref{system}) over $k$.
\end{prop}
\begin{proof}  ({\em cf.,} \cite{CHS07}, Proposition 2.7) If $K$ is a total $\SDP$-PV ring for equations (\ref{system}) over $k$, then Lemma \ref{radical} and Corollary \ref{cornonewconstants} imply the conclusion.

We now assume that $K$ satisfies the three conditions above and let $R$ be a $\SDP$-PV extension of $k$ for the equations (\ref{system}) over $k$.  We will show there is a $k$-$\SDP$ embedding of $R$ into $K$.  Assuming this, we may assume, without loss of generality, that $R \subset K$ and $R = k\{U, \frac{1}{\det U}\}_\Pi$ where $U$ is a fundamental solution matrix for the systems (\ref{system}).  We then have that $Z = U\cdot M$ where $M \in \GL_n(k^\SD)$.  Therefore, $K$ is the total ring of quotients of $R$.

We now will produce a $k$-$\SDP$ embedding of $R$ into $K$. Let $X = (X_{i,j})$ be an $n\times n$ matrix of $\Pi$-differential indeterminates and let $R_0:= k\{X,\frac{1}{\det X}\}_\Pi \subset K\{X,\frac{1}{\det X}\}_\Pi$.  We define a $\SDP$-structure on $K\{X,\frac{1}{\det X}\}_\Pi$ by setting $\sigma(X) = A_iX$ for all $\sigma \in \Sigma$ and $\delta_i(X) = B_i X$ for all $\delta \in \Delta$.  This induces a $\SDP$-structure on $R_0$ as well.  Since $R$ is a simple $\SDP$-ring, there is a maximal $\SDP$-ideal $p \subset R_0$ such that $R$ is isomorphic to $R_0/p$ and, without loss of generality, we may assume $R = R_0/p$. Define elements $Y_{i,j} \in K\{X,\frac{1}{\det X}\}_\Pi$ by the formula $(Y_{i,j}) = Z^{-1}(X_{i,j})$.  Note that for each $i,j$, we have   $\sigma(Y_{i,j}) = Y_{i,j} \ \forall \sigma \in \Sigma$ and $\delta(Y_{i,j}) = 0 \ \forall \delta \in \Delta$ and that $K\{X,\frac{1}{\det X}\}_\Pi = K\{Y,\frac{1}{\det Y}\}_\Pi$ where $Y = (Y_{i,j})$.  The ideal $p \subset R_0 \subset k\{Y,\frac{1}{\det Y}\}_\Pi$ generates a $\SDP$-ideal $P = (p)$ in $K\{X,\frac{1}{\det X}\}_\Pi$.  We define $R_1= k^\SD\{Y,\frac{1}{\det Y}\}_\Pi$ and define $p_1 = P \cap S_1$.  Note that $K\{X,\frac{1}{\det X}\}_\Pi$ satisfies the hypotheses of Lemma \ref{ideals1} and $P$ is a $\SDP$-ideal of this latter ring. Lemma \ref{ideals1} therefore implies that $P$ is generated by $p_1$. Let $m$ be a maximal $\Pi$-ideal of $R_1$ such that $p_1 \subset m$.  Since $k^\SD$ is differentially closed, we have $R_1/m \simeq k^\SD$.  We therefore have a $\Pi$-homomorphism $R_1 \rightarrow R_1/m \simeq k^\SD$.  Since $p_1$ generates $P$ this yields a $\SDP$-homomorphism $K\{X,\frac{1}{\det X}\}_\Pi \rightarrow K$ whose kernel contains $P$.  Restricting to $R_0$, we get a $\SDP$-homomorphism $R_0 \rightarrow K$ whose kernel contains $p$.  Since $p$ is a maximal $\SDP$-ideal, this kernel must equal $p$ and therefore yields a $\SDP$-isomorphism of $R$ into $K$. \end{proof}

\subsubsection{$\SDP$-Galois Groups and the Fundamental Theorem}\label{galoissec}  In this section we shall show that the automorphism group of a $\SDP$-PV ring $R$ is a  $\Pi$-differential group defined over $R^\SD$ and  that there is a Galois correspondence between certain subrings of the total ring of quotients of $R$ and differential subgroups of this automorphism group.  We refer to \cite{cassidy1} and the Appendix of \cite{CaSi} for a review of the elementary language of affine differential geometry and linear differential algebraic groups.

\begin{prop}\label{galoisgp} Let $k$ be a $\SDP$-field and let $C= k^\SD$ be $\Pi$-differentially closed. Let $R$ be a $\SDP$-PV extension of $k$ with total ring of quotients $K$.   Then $\Aut_\SDP(K/k)$ may be identified with the group $G(C)$ of $C$-points of a linear $\Pi$-differential algebraic group $G \subset \GL_n$ defined over $C$.
\end{prop}
\begin{proof}
Let $X = (X_{i,j})$ be $n\times n$ matrix of $\Pi$-differential indeterminates.  We shall assume that $R$ is the $\SDP$--PV extension for the system (\ref{system}) and consider $k\{X,\frac{1}{\det X}\}_\Pi$ as a $\SDP$-ring with the structure defined by $\sigma(X) = A_iX$ for all $\sigma \in \Sigma$ and $\delta(X) = B_iX$ for all $\delta \in \Delta$.  We may therefore write $R = k\{X,\frac{1}{\det X}\}_\Pi/q$ where $q$ is a maximal $\SDP$-ideal. Let $Z$ be the image of $X$ in $R$ and define $Y = (Y_{i,j}) \in K\{X, \frac{1}{\det X}\}_\Pi$ by the formula $X = ZY$.  Note that the $\{Y_{i,j}\}$ are again $\Pi$-differential indeterminates and that $\sigma(Y) = Y$ for all $\sigma \in \Sigma$ and $\delta(Y) = 0$ for all $\delta \in \Delta$.  We will show that there is a radical $\Pi$ ideal $I\subset S = C\{Y,\frac{1}{\det Y}\}_\Pi \subset K\{Y, \frac{1}{\det Y}\}_\Pi = K\{X, \frac{1}{\det X}\}_\Pi$ such that $S/I$ is the coordinate ring of a linear $\Pi$-differential algebraic group $G$ and $\Aut^\SDP(K/k) = G(C)$.

Consider the following rings:

\begin{eqnarray*}\label{ringseqn}k\{X,\frac{1}{\det X}\}_\Pi \subset K\{X, \frac{1}{\det X}\}_\Pi = K\{Y, \frac{1}{\det Y}\}_\Pi \supset C\{Y,\frac{1}{\det Y}\}_\Pi
\end{eqnarray*}
Since $K$ is a direct sum of fields of characteristic zero, the
ideal $qK\{X, \frac{1}{\det X}\}_\Pi \subset K\{X, \frac{1}{\det
X}\}_\Pi $ is a radical $\SDP$-ideal ({\rm cf.,} \cite{PuSi2003},
Corollary A.16).  We now consider this ideal as  an ideal in $K\{Y,
\frac{1}{\det Y}\}_\Pi $ and apply Lemma~\ref{ideals1} to conclude
that $qK\{X, \frac{1}{\det X}\}_\Pi $ is generated by $I = qK\{X,
\frac{1}{\det X}\}_\Pi \cap C\{Y, \frac{1}{\det Y}\}_\Pi$. Clearly
$I$ is a radical $\Pi$-ideal of $S= C\{Y, \frac{1}{\det Y}\}_\Pi$.
We shall show that $S/I$ is the coordinate ring of a linear
$\Pi$-differential algebraic group $G$ and that $G(C)$ may be
identified with $\Aut_\SDP(K/k)$.

The group $\Aut_\SDP(K/k)$ can be identified with the set of $c = (c_{i,j}) \in \GL_n(C)$ such that the map $X \mapsto Xc$ leaves the ideal $q$ invariant.  One can easily show that the following statements are equivalent.
\begin{enumerate}
\item[(i)] $c \in \Aut_\SDP(K/k)$,
\item [(ii)] The map $k\{X, \frac{1}{{\det X}}\}_{\Pi} \rightarrow K$
defined by  $X \mapsto Zc$ maps all elements of $q$
to zero.
\item[(iii)] The map $K\{X,\frac{1}{{\det X}}\}_{\Pi}\rightarrow K$
defined by  $X\mapsto Zc$ maps all elements of
$qK\{X, \frac{1}{\det X}\}_\Pi =qK\{Y, \frac{1}{\det Y}\}_\Pi $
to zero.
\item[(iv)] Considering $qK\{X, \frac{1}{\det X}\}_\Pi$ as an
ideal of $qK\{Y, \frac{1}{\det Y}\}_\Pi$, the map
$K\{Y, \frac{1}{\det Y}\}_\Pi \rightarrow K$ defined by $Y \mapsto c$ sends all
elements of $qK\{Y, \frac{1}{\det Y}\}_\Pi$ to zero.
\end{enumerate}
 Since the ideal
$qK\{Y, \frac{1}{\det Y}\}_\Pi$ is generated by
$I$, the last statement above is equivalent to $c$ being a zero of the
ideal $I$.  Since $\Aut_\SDP(K/k)$ is a group, the set
$G(C)$ is a subgroup of $\GL_n(C)$.  Therefore $G$ is a linear differential algebraic group. \end{proof}

We note that the differential group structure on $\Aut_\SDP(R/k)$ can be defined without assuming that $k^\SD$ is differentially closed.  We can replace the assumption that $k$ is differentially closed with the assumption that $R^\SD = k^\SD = C$.  We then consider the functor $\calG$ from $\Pi$-differential $C$-algebras $B$ to groups given by $\calG(B) = \Aut_\SDP(R\otimes_CB/k\otimes_CB)$ where the $\SDP$ structure on $R\otimes_CB$ and $k\otimes_CB$ are given by $\sigma(f\otimes b) = \sigma (f) \otimes b \ \forall \sigma \in \Sigma$ and $\delta(f\otimes b) = \delta(f) \otimes b \ \forall \delta \in \Delta$ and $\der(f\otimes b) = \der(f)\otimes b + f\otimes \der(b) \ \forall \der \in \Pi$.  One can show (as in \cite{PuSi2003}, Theorem 1.27) that this functor is representable and so defines a linear differential algebraic group ({\rm cf.,} \cite{ovchinnikov1}).  One can furthermore show that when $k^\SD$ is differentially closed,  the linear differential algebraic group structure is the same as that defined above.  Yet another approach to showing that the Galois group is a linear differential algebraic group is to develop a theory of $\SDP$-modules analogous to differential modules or difference modules (as in \cite{PuSi2003} and \cite{PuSi}) and show that the category of such objects forms a differential tannakian category in the sense of \cite{ovchinnikov2}, where it is shown that the automorphism group of the associated differential fibre functor has the structure of a linear differential proalgebraic group.

We now show a weak version of the Fundamental Theorem (see Theorem~\ref{fundthm} below).

\begin{lem}\label{weakfundthm} Let $k$ be a $\SDP$-field with $k^\SD$ $\Pi$-differentially closed. Let $R$ be a  $\SDP$-PV ring for for  the equations (\ref{system}) over $k$ and let  $G = \Aut_\SDP(R/k)$ be the $\SDP$-Galois group.  Let $K$ be the total quotient ring of $R$ and $H$ a Kolchin $\Pi$-closed subgroup of $G$.
\begin{enumerate}
\item The ring $K^G$ of elements of $K$ left fixed by $G$ is $k$.
\item If $K^H = k$ then $H = G$.
\end{enumerate}
\end{lem}
\begin{proof} 1. ({\em cf.,} \cite{PuSi2003}, Theorem 1.27(3))  Let $a = \frac{b}{c} \in K\backslash k$ with $b, c \in R$ and let $d = b\otimes c -c\otimes b$.  Using the description of $R$ given in Lemma~\ref{idempotents} and the fact that $k$ has characteristic $0$, we see that the ring $R \otimes_k R$ has  no zero divisors (\cite{PuSi2003}, Lemma A.16).  Since $d \neq 0$,  we can localize and consider a maximal $\SDP$-ideal  $J$ in $ R \otimes_k R[\frac{1}{d}]$ (with the obvious $\SDP$-structure). Let $S =  R \otimes_k R[\frac{1}{d}]/J$. Since $k^\SD$ is $\Pi$-differentially closed we have $k^\SD = S^\SD$.  Furthermore, the obvious maps $\phi_i:R\rightarrow S$ are injections.  The images $\phi_1(R)$ and $\phi_2(R)$ are generated by fundamental matrices of the same $\SD$-system and so generate the same ring $S_0 \subset S$.  Therefore $\sigma = \phi_2^{-1}\circ \phi_1$ is in $\Aut_\SDP(R/k)$.  The image of $d$ in $S$ is $\phi_1(b)\phi_2(c) - \phi_1(c)\phi_2(b)$ which is nonzero.  Therefore $\phi_1(b)\phi_2(c) \neq \phi_1(c)\phi_2(b)$ and so $\sigma(b)c - \sigma(c)b \neq 0$.  This implies $\sigma(\frac{b}{c}) \neq \frac{b}{c}$.

2. ({\em cf.,} \cite{CaSi}, Proposition 9.10)
 Assuming that $H \neq G$, we shall derive a contradiction.  Let $ R = k\{X, \frac{1}{\det X}\}_\Pi/ p$ where  $X$ is an $n \times n$ matrix of $\Pi$-differential indeterminates and denote by $Z$ the image of $X$ in $R$. We consider $k\{X, \frac{1}{\det X}\}_\Pi \subset K\{X, \frac{1}{\det X}\}_\Pi$ and let $Y = Z^{-1}X$.  Note that $\sigma(Y) = Y$ for all $\sigma \in \Sigma$ and $\delta Y = 0$ for all $\delta \in \Delta$. Let $I \subset C\{Y, \frac{1}{Y}\}_\Pi \subset K\{Y, \frac{1}{\det Y}\}_\Pi =K\{X, \frac{1}{\det X}\}_\Pi$ be the defining ideal of $G$ (where $C = k^\SD$).   Under our assumption that $H \neq G$, we have that there is an element $P \in C\{Y, \frac{1}{Y}\}_\Pi$ such that $P \notin I$ and $P(h) = 0 $ for all $h \in H(C)$. Lemma~\ref{ideals1} implies that $P \notin (I) = IK\{Y, \frac{1}{\det Y}\}_\Pi$.  Let $T = \{Q \in K\{X, \frac{1}{\det X}\}_\Pi \ | \ Q \notin (I) \mbox{ and } Q(Zh) = 0 \mbox{ for all } h \in H\}$.  Note that $P(Z^{-1}X) \in T$, so $T \neq \emptyset$.  Any element of $ K\{Y, \frac{1}{\det Y}\}_\Pi$ can be written as $\sum_\alpha f_\alpha Q_\alpha$ with $f_\alpha\in K$ and $Q_\alpha \in k\{X, \frac{1}{\det X}\}_\Pi$.  Select $Q = f_{\alpha_1}Q_{\alpha_1} + \ldots +f_{\alpha_m}Q_{\alpha_m} \in T$ with $m$ minimal.

 We shall first show that we can there is such a $\tilde{Q} = \tilde{f}_{\alpha_1}\tilde{Q}_{\alpha_1} + \ldots +\tilde{f}_{\alpha_m}\tilde{Q}_{\alpha_m} \in T$ with $m$ minimal and each nonzero $\tilde{f}_{\alpha_i}$ invertible in $K$. To see this, let $Q$ be as above and  note that each $e_iQ$ satisfies $e_iQ(Zh) = 0 \mbox{ for all } h \in H$.  Therefore for some $i$ we must have $e_iQ \notin (I)$ since $Q = \sum e_iQ$.  We will assume $e_0Q \notin (I)$ and therefore have that $e_0Q \in T$ and is minimal.  The elements of $H$ permute the elements $\{e_i\}$ and since $K^H = k$, we must have that $H$ acts transitively on this set.  Therefore the exist $h_i\in H$ such that $h_i(e_0) = e_i$.  Let
 $\tilde{Q} = \sum h_i(e_0Q)$.  We again have that $\tilde{Q} \notin (I)$ (otherwise $e_0Q = e_0\tilde{Q} \in (I)$) and so we can write $\tilde{Q} = \tilde{f}_{\alpha_1}\tilde{Q}_{\alpha_1} + \ldots +\tilde{f}_{\alpha_m}\tilde{Q}_{\alpha_m} \in T$ with $m$ minimal and each nonzero $\tilde{f}_{\alpha_i}$ invertible in $K$.

 Since all the $f_{\alpha_i}$ are invertible,  we may assume that $f_{\alpha_1} = 1$.  For each  $h \in H$, let $Q^h = f^h_{\alpha_1}Q_{\alpha_1} + \ldots +f^h_{\alpha_m}Q_{\alpha_m}$.  Note that $Q^h \in T$ as well.  Since $Q-Q^h$ is shorter than $Q$ and $(Q-Q^h)(Zh) = 0$ for all $h \in H$, we have that $Q-Q^h \in (I)$.  Assume that  $Q - Q^h \neq 0$.

 Let $\calY$ be the monomials appearing with nonzero coefficients in $Q$. The monomials in $Q-Q^h$ will be a proper subset of $\calY$. Using the fact that $H$ acts transitively on the $e_i$ we can find some non zero $S \in (I)$ with $S(Zh) = 0$ for all $h \in H$ as above such that the monomials of $S$ are in $\calY$ and the coefficients of these monomials are invertible in $K$.   There then exists an $l \in K$ such that $Q - lS$ is shorter than $Q$ and must lie in $T$, a contradiction. Therefore $Q = Q^h$ for all $h \in H$ and so $Q \in k\{X,\frac{1}{\det X}\}_\Pi$.  Since $Q(Z\cdot id) = 0$ we have $Q(Z\cdot g) = g(Q(Z\cdot id) = 0$ for all $g\in G$, contradicting the fact that $Q \notin (I)$. \end{proof}

As noted in Section 1.3 of \cite{PuSi}, one cannot have a Galois correspondence between subrings of a Picard-Vessiot extension of a difference equation and closed subgroups of the Galois group and one must look at subrings of the total ring quotients. In the next proposition, we characterize certain $\SDP$-subrings of the total ring of quotients of a $\SDP$-PV ring that will be in bijective correspondence with the $\Pi$-differentially closed subgroups of $\Aut_\SDP(R/k)$.
\begin{thm} \label{fundthm} Let $k$ be a $\SDP$- field with $C = k^\SD$ differentially closed.  Let $K$ be a total $SDP$-PV ring over $k$ and $G= \Aut_\SDP(K/k)$.  Let $\calF$ be the set of $\SDP$-rings $F$ with $k \subset F \subset K$ such that every non zero divisor of $F$ is a unit in $F$.  Let $\calG$   denote the set of $\Pi$-differential algebraic subgroups of $G$.  \begin{enumerate}
\item For any $F \in \calF$, the subgroup $G(K/F) \subset G$ of elements of $G$ which fix $F$ pointwise is a $\Pi$-differential subgroup of $G$.
\item For any differential algebraic subgroup $H \subset G$, the ring $K^H$ of elements left fixed by $H$ belongs to $\calF$.
\item Let $\alpha:\calF \rightarrow \calG$ and $\beta:\calG \rightarrow \calF$ denote the maps $F \mapsto G(K/F)$ and $H \mapsto K^H$.  Then $\alpha$ and $\beta$ are each other's inverses.
\end{enumerate}
\end{thm}

\begin{proof} ({\em cf.,} \cite{PuSi}, Theorem 1.29) Statement 2.~is clear.  To verify statement 1, note that any element of $ f = \frac{g}{h} \in F$ is the quotient of two $\Pi$-differential polynomials $g,h$ in $Z = (z_{i,j})$, the fundamental matrix for the $\SD$-system associated with $K$.  The condition that $\phi  = (\phi_{i,j})\in G$ leaves $f$ fixed is that $g(\phi(Z)) h(Z) - g(Z)h(\phi(Z)) = 0$ and this is equivalent to a system of $\Pi$-differential polynomial equations in the $\phi_{i,j}$ with coefficients in $C$.

To prove Statement 3., it is clear that $F \subset \beta\alpha(F)$ and $H \subset \alpha\beta(H)$.  We will  show equality in both cases.  Let $K$ be the total ring of quotients of a $\SDP$-PV ring $R$ and let $e_i$ and $R_i$ be as in Lemma~\ref{idempotents}.  Corollary~\ref{zerodivisor} implies that for each $F \in \calF$ we may write $F = \oplus_{i=0}^{d-1} FE_i$ where each $FE_i$ is a field invariant under $\sigma^d$ for all $\sigma \in \Sigma$.  One can further more show that:
\begin{itemize}
\item The ring $KE_i$ is the total $\tilde{\Sigma}\Delta\Pi$-PV ring over $k$  where $\tilde{\Sigma} = \{\tilde{\sigma }= \sigma^d \ | \ \sigma \in \Sigma \}$ for the system $\tilde\sigma_i(Y) = \tilde{A}_i Y, \ \tilde{A_i} = \sigma_i^{d-1}(A_i) \ldots \sigma_i^2(A_i) \sigma_i(A_i) A_i$ for all $\tilde{\sigma}_i \in \tilde{\Sigma}$. This follows from the fact that $(KE_i)^{\tilde{\Sigma}\Delta\Pi }= k^{\tilde{\Sigma}\Delta\Pi} = k^\SDP$ and Proposition~\ref{totalquotprop}.
\item The elements of $G(K/F)$ can be described as tuples $(\phi_0, \ldots , \phi_{d-1})$ where each $\phi_j$ is a $\tilde{\Sigma}\Delta\Pi$-automorphism of $KE_i$ over $FE_i$ and $\sigma\phi_j = \phi_{\sigma(j)}\sigma$ for all $\sigma \in \Sigma$
\end{itemize}
Lemma~\ref{weakfundthm}.1 applied to each $FE_i \subset KE_i$ shows that the set of $G(K/F)$-invariant elements are of $K$ are $F$.  In a similar way,  we may write $K^H =  \oplus_{i=0}^{d-1} K^HE_i$ for some $d$ and apply Lemma~\ref{weakfundthm}.2 to to each $K^H E_i\subset KE_i$ and conclude that $H = G(K/K^H)$.
\end{proof}

We end this section with results that compare the $\SDP$-PV theory with the usual Picard-Vessiot theories of differential and difference equations.  Let $k$ be a $\SDP$ field with $C = k^\SD$ differentially closed. Once again consider the system (\ref{system}) above.  Let $R = k\{Z,\frac{1}{Z}\}_\Pi$ be the $\SDP$-PV ring.  The subring $S = k[Z,\frac{1}{\det Z}]$ is an $\SD$- subring and one can ask if this is a $\SD$-PV ring for the the system (\ref{system}) (or more precisely a $\SD\tilde{\Pi}$-PV ring where $\tilde{\Pi} = \emptyset$, but we will not use this latter infelicitous nomenclature). The following proposition answers this affirmatively, generalizes this  and compares to two relevant Galois groups.   We will denote by $\Theta$ the semigroup generated by $\Pi$, that is the set of derivative operators as in \cite{DAAG}, and by $\Theta(s)$ the set of elements of $\Theta$ of order less than or equal to $s$.

\begin{prop}\label{propcompare} Let $k$ be a $\SDP$ field with $C = k^\SD$ differentially closed and let $R = k\{Z,\frac{1}{Z}\}_\Pi$ be the $\SDP$-PV ring for (\ref{system}).  Fix an integer $s$ and let $\Theta(s)Z = \{ \theta(Z) \ | \ \theta \in \Theta(s)\}$. Let $S = k[ \Theta(s)Z, \frac{1}{\det Z}]$.  Then
\begin{enumerate}
\item $S$ is a $\SD$-PV ring for some $\SD$-system, and
\item the $\SDP$-Galois group $\Aut_\SDP(R/k)$ is Zariski dense in $\Aut_\SD(S/k)$.
\end{enumerate}
\end{prop}
\begin{proof} To simplify notation, we will assume that $\Pi = \{\partial\}$ is a singleton.  We note that for any positive integer $i$,
\begin{eqnarray*}
\sigma_\ell(\partial^iZ) = \partial^i(\sigma Z) & =& \sum_{j=0}^i \binom{i}{j} \partial^j(A_\ell)\partial^{i-j}(Z), \mbox{ for all } \sigma_\ell \in \Sigma\\
\delta_\ell(\partial^iZ) = \partial^i(\delta_\ell Z) & = & \sum_{j=0}^i \binom{i}{j} \partial^j(B_\ell)\partial^{i-j}(Z), \mbox{ for all } \delta_\ell \in \Delta
\end{eqnarray*}
Therefore the matrix
\[U = \left(\begin{array}{ccccc} Z & 0 & 0& \cdots & 0\\ \partial Z & Z & 0 & \ldots & 0 \\ \partial^2 Z & \partial Z & Z & \ldots & 0 \\ \vdots & \vdots & \vdots & \vdots & \vdots \\
\partial^{s-1}Z & \partial^{s-2} Z & \partial^{s-3} Z & \ldots & 0 \\
\partial^s Z & \partial Z^{s-1} & \partial^{s-2}Z & \ldots & Z
 \end{array} \right) \]
 satisfies the systems
 \begin{eqnarray*}
\sigma_\ell(U) &=& \overline{A}_\ell U, \ \mbox{ for all } \sigma_\ell \in \Sigma\\
\delta_\ell(U) &= & \overline{B}_\ell U, \ \mbox{ for all } \delta_\ell \in \Delta
\end{eqnarray*}
where
\[\overline{A}_\ell = \left(\begin{array}{ccccc} A_\ell & 0 & 0& \cdots & 0\\ \binom{s}{1}\partial A_\ell & A_\ell & 0 & \ldots & 0 \\ \binom{s}{2}\partial^2 A_\ell &\binom{s}{1} \partial A_\ell & A_\ell & \ldots & 0 \\ \vdots & \vdots & \vdots & \vdots & \vdots \\
\binom{s}{s-1}\partial^{is-1}A_\ell & \binom{s}{s-2}\partial^{s-2} A_\ell & \binom{s}{s-3}\partial^{s-3} A_\ell & \ldots & 0 \\
\binom{s}{s}\partial^s A_\ell & \binom{s}{s-1}\partial A_\ell^{s-1} & \binom{s}{s-2}\partial^{s-2}A_\ell & \ldots & A_\ell
 \end{array} \right)\]
 and $\overline{B}_\ell$ is a similar matrix with the $A_\ell$ replace by $B_\ell$.

 Corollary~\ref{zerodivisor} implies that the total quotient ring $K_S$ of $S$ can be embedded in the total quotient ring $K_R$ of $R$.  Corollary~\ref{nonewconstants} implies that $K_R^\SD = C$ so we have that $K_S^\SD = C$.  Proposition~\ref{totalquotprop} (for an empty $\Pi$) implies that $S$ is therefore the $\SD$-PV ring for system (\ref{system}).  This proves statement 1.

The group $\Aut_\SDP(R/k)$ leaves $S$ invariant and its action on $R$ is determined by its action on $S$, so we may consider   $\Aut_\SDP(R/k)$ as a subgroup of $\Aut_\SD(S/k)$. Both act on $K_S$ and have $k$ as their fixed fields so by Theorem~\ref{fundthm}, $\Aut_\SDP(R/k)$ is Zariski dense in $\Aut_\SD(S/k)$ (the $\Pi$-Kolchin topology when $\Pi$ is empty is the Zariski topology).\end{proof}

\begin{cor} Let $k,R$ be as above. The $R$ is a simple $\SD$-ring. \end{cor}
\begin{proof} Let $I$ be a $\SD$-deal in $R$.  If $I \neq (0)$, then there is an integer $s$ such that $I_s = I \cap k[\Theta(s)Z, \frac{1}{\det Z}] \neq (0)$.  Since $k[\Theta(s)Z, \frac{1}{\det Z}]$ is a $\SD$-PV ring, it is a simple $\SD$ ring, so $1 \in I_S \subset I$. \end{proof}

\subsubsection{Torsors}\label{torsorsec}  In this section we show that a $\SDP$-PV ring $R$ is a $\Pi$-torsor for the $\SDP$-Galois group $\Aut_\SDP(R/k)$. For basic facts concerning differential algebraic sets see \cite{cassidy1, CaSi}.

\begin{defin} Let $k$ be a $\Pi$-field and $G$ a linear differential algebraic group
defined over $k$. A {\em $G$-torsor (defined over $k$)} is an  affine differential algebraic variety
$V$ defined over $k$ together with  a differential polynomial map $f:V\times_k G \rightarrow
V\times_k
V$ (denoted by  $f:(v,g) \mapsto vg$) such that
\begin{enumerate}
\item  for any $\Pi$-$k$-algebra $S\supset k, v\in V(S), g,g_1,g_2 \in G(S)$, $v1_G = v$, $v(g_1g_2) =
(vg_1)g_2$ and
\item  the
associated homomorphism $k\{V\}\otimes_k k\{V\} \rightarrow k\{V\}\otimes_k k\{G\}$ is an
isomorphism (or equivalently, for any $S \supset k$, the map $V(S)\times G(S) \rightarrow
V(S)\times V(S)$ is a bijection.
\end{enumerate}\end{defin}

For $G \subset \GL_n$, we shall be interested in $G$-torsors  $V$ that come with a natural inclusion $V \subset \GL_n$ and where the action of $G$ on $V$ is given by multiplication on the right.  Assume this is the case and let $S$ is a $\Pi$-$k$-algebra such that $V(S) \neq \emptyset$, then, just as in the situation of torsors over an algebraic group, $V$ is a $G$-torsor if and only if $V(S) = pG(S)$ for some $p \in \GL_n(S)$ ({\em cf.}, \cite{PuSi2003}, \cite{waterhouse}).

\begin{prop}\label{torsorprop} ({\em cf}., \cite{PuSi2003}, Theorem 1.28) Let $k$ be a $\SDP$-field with $C = k^\SD$ $\Pi$-differentially closed.  Let $R$ be a $\SDP$-PV extension of $k$ and $G$ the $\SDP$-Galois group of $R$ over $k$.  Then $R$ is the coordinate ring of a $G$-torsor over $k$.
\end{prop}
\begin{proof}  We will use the notation of Proposition~\ref{galoisgp} and its proof.  In particular, we let $K$ be the total ring of quotients of $R$, $X$ an $n\times n$ matrix of $\Pi$-differential variables, $k\{X, \frac{1}{\det X}\}_\Pi$ a $\SDP$ ring with structure defined as in that proposition $R = k\{X, \frac{1}{\det X}\}_\Pi/q$, $q$ a maximal $\SDP$-ideal and $Y =Z^{-1}X$. In the sequence of rings

\begin{eqnarray*}\label{ringseqn}k\{X,\frac{1}{\det X}\}_\Pi \subset K\{X, \frac{1}{\det X}\}_\Pi = K\{Y, \frac{1}{\det Y}\}_\Pi \supset C\{Y,\frac{1}{\det Y}\}_\Pi
\end{eqnarray*}
we have that $\Sigma$ and $\Delta$ act trivially on $C\{Y,\frac{1}{\det Y}\}_\Pi$.  The action of $\Aut_\SDP(K/k)$ on $K$ gives an action of this group on $K\{X, \frac{1}{\det X}\}_\Pi$ by letting the group act trivially on $X$.  Note that for each $\phi \in \Aut_\SDP(K/k)$ there is a $M_\phi \in \GL_n(C)$ such that $\phi(Z) = Z M_\phi$.  We therefore have that $\phi(Y) = M_\phi^{-1}Y$ and so have an action of $\Aut_\SDP(K/k)$ on $C\{Y,\frac{1}{\det Y}\}_\Pi$.  In Lemma~\ref{ideals2} below we shall show that there is a bijection $I \mapsto (I)$ from the set of $\SDP$-ideals of $ k\{X, \frac{1}{\det X}\}_\Pi$ to the $\Aut_\SDP(K/k)$-invariant $\SDP$-ideals of $K\{X, \frac{1}{\det X}\}_\Pi$.  From Lemma~\ref{ideals1}, we know that the map $J\mapsto (J)$ is a bijection from the set of  $\Aut_\SDP(K/k)$-invariant $\Pi$-ideals of $C\{Y,\frac{1}{\det Y}\}_\Pi$ to the $\Aut_\SDP(K/k)$-invariant $\SDP$-ideals of $K\{Y, \frac{1}{\det Y}\}_\Pi = K\{X, \frac{1}{\det X}\}_\Pi$.  Therefore we have a bijection between the $\SDP$-ideals of $ k\{X, \frac{1}{\det X}\}_\Pi$ and the $\Aut_\SDP(K/k)$-invariant $\Pi$-ideals of $C\{Y,\frac{1}{\det Y}\}_\Pi$. The maximal $\SDP$-ideal $q\subset k\{X,\frac{1}{\det X}\}_\Pi $ lifts to a maximal $\Aut_\SDP(K/k)$ invariant $\SDP$-ideal  $qK\{X,\frac{1}{\det X}\}_\Pi $ which then restricts to  a maximal $\Aut_\SDP(K/k)$ invariant $\Pi$-ideal $r:= qK\{X,\frac{1}{\det X}\}_\Pi \cap C\{Y,\frac{1}{\det Y}\}_\Pi$.  By maximality, $r$ is a radical ideal so its zero set $W$ is a minimal $\Aut_\SDP(K/k)$ invariant $\Pi$-Kolchin closed subset of $\GL_n(C)$. It therefore must be a coset of $\Aut_\SDP(K/k)$. Since $X=Z$ is a zero of $q$, we have that $Y = id$, where $id$ is the identity matrix, is a zero of $r$, that is $id \in W$.  We therefore have that $W=\Aut_\SDP(K/k)$.  In particlar, this implies that if $V$ is the $\Pi$-closed subset of $\GL_n$ defined by $q$, then $V(K) = ZG(K)$ where $G$ is the linear differnetial algebraic group whose $C$-points $G(C)$ are $\Aut_\SDP(K/k)$.\end{proof}
\begin{lem}\label{ideals2} The map $I \mapsto (I)$ from the set of  ideals in $k\{X, \frac{1}{\det X}\}_\Pi$ to the set of $\Aut_\SDP(K/k)$-invariant ideals of $K\{X, \frac{1}{\det X}\}_\Pi$ is a bijection.\end{lem}
\begin{proof}  We modify Lemma 1.29 of \cite{PuSi2003} in the same way that Lemma~\ref{ideals1} was a modification of Lemma 1.23 of \cite{PuSi2003}.  We wish to show that any $\Aut_\SDP(K/k)$-invariant ideal $J$ of $K\{X, \frac{1}{\det X}\}_\Pi$ is generated by $I:=J\cap k\{X, \frac{1}{\det X}\}_\Pi$.  Let $\{b_\alpha\}$ be a $k$-basis of $k\{X, \frac{1}{\det X}\}_\Pi$ and, for any $f \in J$, let $f = \sum f_\alpha b_\alpha$, with the $f_\alpha \in K$.  As before, define the length  $\ell(f)$ of $f$ to be the number of $\alpha$ such that $f_\alpha \neq 0$.  We will show by induction on the length of $f$ that $f \in (I)$.  We can assume that $\ell(f) = \ell >1$ and that the claim is true for all elements of length smaller than $\ell(f)$.

Let us  assume that we can show that $f \in (I)$ if $f$has the property  $f = \sum_{i=1}^{t-1} \phi_i(e_j h$)  where $h \in K$ and $\phi_i \in \Aut_\SDP(K/k)$ and satisfy $\phi_i(e_j) = e_i,i \neq j, \ \phi_j = id$.   For a general $g \in J$, to show that $g \in (I)$ it is enough to show that each $e_jg \in (I)$.  We know that $\Aut_\SDP(K/k)$ permutes the $e_i$ since   and acts transitively on these since the fixed field is $k$.  Select $\phi \in \Aut_\SDP(K/k)$ such that $\phi_j = id$ and $\phi_i(e_j) = e_i$ for $i \neq j$.  Let $f = \sum \phi_i(e_jg)$.  Assuming we can show $f \in (I)$, we have $e_jf = e_jg \in I$.

Therefore we shall assume that $f$ is of the above form.  Note
that  if $f_\alpha \neq 0$, then $f_\alpha$ is invertible in $K$ and so we can assume that $f_{\alpha_1} = 1$ for some $f_{\alpha_1}$.  If all $f_\alpha \in k$ we are done so assume there exists a $f_{\alpha_2} \in K\backslash k$.  For any $\phi \in \Aut_\SDP(K/k)$, we have that $\ell(\phi(f) - f) < \ell(f)$ so $\phi(f) - f \in (I)$.  Since the fixed field of  $\Aut_\SDP(K/k)$ is $k$, we have that there is some $\phi \in \Aut_\SDP(K/k)$ such that  $\phi(f_{\alpha_2}) \neq f_{\alpha_2}$.  One sees that $\phi(f_{\alpha_2}^{-1}f) - f_{\alpha_2}^{-1}f \in (I)$ and so $(\phi(f_{\alpha_2}^{-1} )- f_{\alpha_2}^{-1})f =  (\phi(f_{\alpha_2}^{-1}f) - f_{\alpha_2}^{-1}f)  - \phi(f_{\alpha_2}^{-1})(\phi(f)-f ) \in (I)$.  Since $(\phi(f_{\alpha_2}^{-1} )- f_{\alpha_2}^{-1})$ may not be invertible, we may not yet conclude that $f \in (I)$ but we can conclude that for some $i$, $e_i f \in (I)$.  We then have that for any $l, 0\leq l \leq t-1 \  \phi_t(\phi_i^{-1}(e_if)) = \phi_t(\phi_i^{-1}(\phi_i(e_jh))) = \phi_t(e_jh) = e_tf \in (I)$ and so $f \in (I)$.
\end{proof}

An immediate consequence of this last result concerns the differential dimension of a $\SDP$-PV extension.  Let $k$ be a $\Pi$-field  (of characteristic 0) and $R$ a $k$-$\Pi$ algebra.  We say elements $r_1, \ldots, r_m \in R$ are  {\em$\Pi$-differentially dependent over $k$} if there exists a nonzero differential polynomial $P\in k\{Y_1, \ldots , Y_m\}_\Pi$ such that $P(r_1, \ldots , r_n) = 0$.  Elements that are not differentially dependent are said to be {\em differentially independent}.  If $K$ is $\Pi$-differential field, the {\em $\Pi$-differential transcendence degree of $K$ over $k$} $\Pi$-diff.tr.deg.$(K/k)$ is the size of a  maximal differentially independent subset of $K$.  If $p$ is a prime $\Pi$-ideal in $k\{Y_1, \ldots , Y_m\}_\Pi$, the {\em $\Pi$-dimension of $p$ over k} is defined to be $\Pi$-diff.tr.deg.$(K/k)$ where $K$ is the quotient field of $k\{Y_1, \ldots , Y_m\}_\Pi/p$.  If $q$ is a radical $\Pi$-ideal, the $\Pi$-dimension of $q$ over $k$ is the maximum of the $\Pi$-dimensions over $k$ of the prime differential  ideals $p_i$ where $q = \cap p_i$. If $R$ is  a finitely generated reduced $\Pi$-$k$-algebra, the the {\em $\Pi$-dimension of $R$ over $k$, $\Pi$-dim.$_k(R)$},  is defined to be the $\Pi$-dimension over $k$ of the ideal $q$ where $R = k\{Y_1, \ldots , Y_m\}_\Pi/q$ and if $V$ is a Kolchin closed set, defined over $k$, the  differential dimension of $V$ over $k$ is the $\Pi$-dimension of a radical differential ideal defining it.  In \cite{DAAG}, Kolchin shows (using differential dimension polynomials) that if $k \subset k'$ are differential fields and $q$ is a radical differential ideal in $k\{Y_1, \ldots , Y_m\}_\Pi$, then the differential dimension of $q$  (over $k$) and $qk'\{Y_1, \ldots , Y_m\}_\Pi$ (over $k'$) are the same.  The following result easily now follows from Proposition~\ref{torsorprop}.

\begin{prop}\label{dimprop} Let $k$ be a $\SDP$-field with $C = k^\SD$ $\Pi$-differentially closed.  Let $R$ be a $\SDP$-PV extension of $k$ and $G$ the $\SDP$-Galois group of $R$ over $k$. Then $\Pi$-dim$._k(R)$ = $\Pi$-dim.$_C(C\{G\}_\Pi)$.\end{prop}
\begin{proof} Proposition~\ref{torsorprop} states that $R$ is the differential coordinate ring of a $G$-torsor  $V$ over $k$.  Let $\tilde{k}$ be a $\Pi$-differentially closed field containing $k$.  Since $V(\tilde{k}) \neq \emptyset$ we have that $V(\tilde{k}) \simeq G(\tilde{k})$, that is $\tilde{k}\otimes_k R \simeq \tilde{k} \otimes_C k\{G\}_\Pi$.  The result now follows from the discussion preceding this proposition.\end{proof}

We note that Proposition~\ref{torsorprop} can be used to give another proof of Theorem~\ref{fundthm} as in Theorem 1.29 of \cite{PuSi}.  To generalize the proof of this latter result one needs to know that if $H$ is a proper differential subgroup of a linear differential algebraic group $G$, then there exists a nonconstant differential rational function on $G$ that is left invariant by the action of $H$ coming from right multiplication ({\em cf.,} Proposition 14 of \cite{cassidy1}).

\addcontentsline{toc}{section}{References}

\begin{thebibliography}{10}

\bibitem{abramov75}
S.~A. Abramov.
\newblock The rational component of the solution of a first order linear
  recurrence relation with rational right hand side.
\newblock {\em \v Z. Vy\v cisl. Mat. i Mat. Fiz.}, 15(4):1035--1039, 1090,
  1975.

\bibitem{abramov_zima}
S.A. Abramov and E.~V. Zima.
\newblock {D'Alembert} solutions of inhomogeneous linear equations
  (differential, difference and otherwise).
\newblock In {\em Proceedings of the 1996 International Symposium on Symbolic
  and Algebraic Computation}, pages 232--239. ACM Press, 1996.

\bibitem{andre_galois}
Y. ~Andr{\'e}.
\newblock Diff\'erentielles non commutatives et th\'eorie de {G}alois
  diff\'erentielle ou aux diff\'erences.
\newblock {\em Ann. Sci. \'Ecole Norm. Sup. (4)}, 34(5):685--739, 2001.

\bibitem{bankkaufman}
S.~Bank and R.~Kaufman.
\newblock A note on {H\"older's} theorem concerning the gamma function.
\newblock {\em Math. Ann.}, 232:115--120, 1978.

\bibitem{barkatou99}
M.~A. Barkatou.
\newblock On rational solutions of systems of linear differential equations.
\newblock {\em Journal of Symbolic Computation}, 28(4/5):547--568, 1999.

\bibitem{bronstein92}
M.~Bronstein.
\newblock On solutions of linear ordinary differential equations in their
  coefficient field.
\newblock {\em Journal of Symbolic Computation}, 13(4):413 -- 440, 1992.

\bibitem{bron_diff}
M.~Bronstein.
\newblock On solutions of linear ordinary difference equations in their
  coefficient field.
\newblock {\em J. Symbolic Comput.}, 29(6):841--877, 2000.

\bibitem{blw}
M.~Bronstein, Z.~Li, and M.~Wu.
\newblock {Picard-Vessiot extensions for linear functional systems}.
\newblock In Manuel Kauers, editor, {\em Proceedings of the 2005 International
  Symposium on Symbolic and Algebraic Computation (ISSAC 2005)}, pages 68--75.
  ACM Press, 2005.

\bibitem{cassidy1}
P.J.~Cassidy.
\newblock Differential algebraic groups.
\newblock {\em Amer. J. Math.}, 94:891--954, 1972.

\bibitem{cassidy2}
P.J.~Cassidy.
\newblock The differential rational representation algebra on a linear
  differential algebraic group.
\newblock {\em J. Algebra}, 37(2):223--238, 1975.

\bibitem{cassidy6}
P.J.~Cassidy.
\newblock The classification of the semisimple differential algebraic groups
  and the linear semisimple differential algebraic {L}ie algebras.
\newblock {\em J. Algebra}, 121(1):169--238, 1989.

\bibitem{CaSi}
P.J.~Cassidy and M.~F. Singer.
\newblock Galois theory of parameterized differential equations and linear
  differential algebraic groups.
\newblock In D.~Bertrand, B.~Enriquez, C.~Mitschi, C.~Sabbah, and R.~Schaefke,
  editors, {\em Differential Equations and Quantum Groups}, volume~9 of {\em
  IRMA Lectures in Mathematics and Theoretical Physics}, pages 113-- 157. EMS
  Publishing House, 2006.

\bibitem{carmichael}
R.D. Carmichael.
\newblock On transcendentally transcendental functions.
\newblock {\em Trans. A.M.S.}, 14(3):311--319, 1913.

\bibitem{CHS07}
Z.~Chatzidakis, C.~Hardouin, and M.~F. Singer.
\newblock {On the Definition of Difference Galois Groups}.
\newblock To appear in the Proceedings of the Newton Institute 2005 semester
  `Model theory and applications to algebra and analysis' ; also {\tt
  arXiv:math.CA/0705.2975}, 2007.

\bibitem{etingof}
P.~I. Etingof.
\newblock Galois groups and connection matrices of {$q$}-difference equations.
\newblock {\em Electron. Res. Announc. Amer. Math. Soc.}, 1(1):1--9
  (electronic), 1995.

\bibitem{hardouin06}
C.~Hardouin.
\newblock Hypertranscendance et groupes de Galois aux diff\'erences.
\newblock {\tt arXiv:math.RT/0702846v1}, 2006.

\bibitem{hardouin07}
C.~Hardouin.
\newblock Hypertranscendance des syst\`emes diagonaux aux diff\'erences.
\newblock to appear in Compositio Mathematica, 2007.

\bibitem{HaSi}
C.~Hardouin and M.~F. Singer.
\newblock Differential independence of solutions of a class of q-hypergeometric
  difference equations.
\newblock Maple worksheet available at  {\tt http://www4.ncsu.edu/$\sim$singer/ms\_papers.html}, 2007.

\bibitem{hausdorff}
F.~Hausdorff.
\newblock {Zum H\"olderschen Satz \"uber $\Gamma(x)$ }.
\newblock {\em Math. Ann.}, 94:244--247, 1925.

\bibitem{hendriks_qdiff}
P.~A. Hendriks.
\newblock An algorithm for computing the standard form for second order linear
  $q$-difference equations.
\newblock {\em Journal of Pure and Applied Mathematics}, 117-118:331--352,
  1997.

\bibitem{hoeij98b}
M.~{\SortNoop{Hoeij}}van~Hoeij.
\newblock Rational solutions of linear difference equations.
\newblock In O.~Gloor, editor, {\em Proceedings of ISSAC'98}, pages 120--123.
  ACM Press, 1998.

  \bibitem{hoeij_sing}
M.~{\SortNoop{Hoeij}}van~Hoeij.
\newblock Finite singularities and hypergeometric solutions of linear
  recurrence equations.
\newblock {\em J. Pure Appl. Algebra}, 139(1-3):109--131, 1999.
\newblock Effective methods in algebraic geometry (Saint-Malo, 1998).



\bibitem{hoelder}
O~H\"older.
\newblock {\"Uber die Eigenschaft der Gamma Funktion keiner algebraische
  Differentialgleichung zu gen\"ugen}.
\newblock {\em Math. Ann.}, 28:248--251, 1887.

\bibitem{ishizaki}
K.~Ishizaki.
\newblock Hypertranscendency of meromorphic solutions of linear functional
  equation.
\newblock {\em Aequationes Math.}, 56:271--283, 1998.

\bibitem{kaplansky}
I.~Kaplansky.
\newblock {\em {An Introduction to Differential Algebra}}.
\newblock Hermann, Paris, second edition, 1976.

\bibitem{karr}
M.~Karr.
\newblock Theory of summation in finite terms.
\newblock {\em J. Symbolic Comput.}, 1(3):303--315, 1985.

\bibitem{kolchin_ostrowski}
E.~R. Kolchin.
\newblock Algebraic groups and algebraic dependence.
\newblock {\em American Journal of Mathematics}, 90:1151--1164, 1968.

\bibitem{DAAG}
E.~R. Kolchin.
\newblock {\em {Differential Algebra and Algebraic Groups}}.
\newblock Academic Press, New York, 1976.

\bibitem{matusevich}
L.~F.~Matusevich.
\newblock Rational summation of rational functions.
\newblock {\em Beitr\"age Algebra Geom.}, 41(2):531--536, 2000.

\bibitem{moore}
E.H. Moore.
\newblock Concerning transcendentally transcendental functions.
\newblock {\em Math. Ann.}, 48:49--74, 1897.

\bibitem{ostrowski25}
A.~Ostrowski.
\newblock {Zum H\"olderschen Satz \"uber $\Gamma(x)$ }.
\newblock {\em Math. Ann.}, 94:1--13, 1925.

\bibitem{ovchinnikov1}
A.~Ovchinnikov.
\newblock Tannakian approach to linear differential algebraic groups.
\newblock {\tt arXiv:math.RT/0702846v1}, 2007.

\bibitem{ovchinnikov2}
A.~Ovchinnikov.
\newblock Tannakian categories, linear differential algebraic groups, and
  parameterized linear differential equations.
\newblock {\tt arXiv:math.RT/0703422v1}, 2007.

\bibitem{paule2}
P.~Paule.
\newblock Greatest factorial factorization and symbolic summation.
\newblock {\em J. Symbolic Comput.}, 20(3):235--268, 1995.

\bibitem{PWZ}
M.~Petkovsek, H.~Wilf, and D.~Zeilberger.
\newblock {\em {A=B}}.
\newblock A. K. Peters, Wellsely, Massachusets, 1996.

\bibitem{praagman}
C.~Praagman.
\newblock Fundamental solutions for meromorphic linear difference equations in
  the complex plane, and related problems.
\newblock {\em J. Reine Angew. Math.}, 369:101--109, 1986.

\bibitem{PuSi}
M. {\SortNoop{Put}}van~der Put and M.~F. Singer.
\newblock {\em {Galois Theory of Difference Equations}}, volume 1666 of {\em
  Lecture Notes in Mathematics}.
\newblock Springer-Verlag, Heidelberg, 1997.

\bibitem{PuSi2003}
M. {\SortNoop{Put}}van~der Put and M.~F. Singer.
\newblock {\em {Galois Theory of Linear Differential Equations}}, volume 328 of
  {\em Grundlehren der mathematischen Wissenshaften}.
\newblock Springer, Heidelberg, 2003.

\bibitem{roques}
J.~Roques.
\newblock Galois groups of basic hypergeometric equations.
\newblock preprint, {\tt arXiv:math.CA/0709.3275v1}, 2007.

\bibitem{rosenlicht_rank}
M.~Rosenlicht.
\newblock The rank of a {H}ardy field.
\newblock {\em Trans. Amer. Math. Soc.}, 280(2):659--671, 1983.

\bibitem{rubel_trans}
L.A. Rubel.
\newblock A survey of transcendentally transcendental functions.
\newblock {\em Am. Math. Monthly}, 96(9):777--788, 1989.

\bibitem{schneider_phd}
C.~Schneider.
\newblock {\em {Symbolic Summation in Difference Fields}}.
\newblock PhD thesis, RISC, J. Kepler University Linz, May 2001.
\newblock (published as Technical report no. 01-17 in RISC Report Series.).

\bibitem{schneider_bounds}
C.~Schneider.
\newblock A collection of denominator bounds to solve parameterized linear
  difference equations in {$\Pi\Sigma$}-extensions.
\newblock {\em An. Univ. Timi\c soara Ser. Mat.-Inform.}, 42(Special issue
  2):163--179, 2004.

\bibitem{schneider_telescope}
C.~Schneider.
\newblock {Parameterized telescoping proves algebraic independence of sums}.
\newblock In {\em {Proceedings of the 19th international conference on formal
  power series and algebraic combinatorics, FPSAC'07}}, pages 1--12, 2007.

\bibitem{seidenberg52}
A.~Seidenberg.
\newblock Some basic theorems in differential algebra (characteristic {$p$},
  arbitrary).
\newblock {\em Trans. Amer. Math. Soc.}, 73:174--190, 1952.

\bibitem{sei58}
A.~Seidenberg.
\newblock Abstract differential algebra and the analytic case.
\newblock {\em Proc. Amer. Math. Soc.}, 9:159--164, 1958.

\bibitem{sei69}
A.Seidenberg.
\newblock Abstract differential algebra and the analytic case. {II}.
\newblock {\em Proc. Amer. Math. Soc.}, 23:689--691, 1969.

\bibitem{um06}
H.~Umemura.
\newblock Galois theory and Painlev\'e equations.
\newblock In \'E. Delabaere and M. Loday-Richaud, eds.,{\em Th\'eories Asymptotiques et \'Equations de Painlev\'e},  Num\'ero 14, S\'eminaires \& Congr\`es, p. 299-340, Soci\'et\'e Math\'ematique de France, 2006.


\bibitem{waterhouse}
W.~C. Waterhouse.
\newblock {\em {Introduction to Affine Group Schemes}}, volume~66 of {\em
  Graduate Texts in Mathematics}.
\newblock Springer-Verlag, New York, 1979.

\bibitem{wu05}
M.~Wu.
\newblock {\em On Solutions of Linear Functional Systems and Factorization of
  Modules over Laurent-Ore Algebras}.
\newblock Ph.d. thesis, {Chinese Academy of Sciences and l'Universit\'e de
  Nice-Sophia Antipolis}, 2005.

\bibitem{zariski_samuel}
O.~Zariski and P.~Samuel.
\newblock {\em {Commutative Algebra}}, volume~1.
\newblock D. Nostrand, Princeton, 1956.

\end{thebibliography}
\newcommand{\SortNoop}[1]{}\def\cprime{$'$}

\end{document}